%% file: agt-2-12.tex
\documentclass{gtart}

\input agtout

\lognumber{12}
\volumenumber{2}
\volumeyear{2002}
\papernumber{12}
\published{16 April 2002}
\pagenumbers{257}{284}
\received{23 July 2001}
\revised{3 April 2002}
\accepted{4 April 2002}

\usepackage{amssymb} 
\usepackage{amsmath} 

\newtheorem{thm}{Theorem}[section]    % Standard theorem environment
\newtheorem{lem}[thm]{Lemma}          % Lemma environment with numbering 
\newtheorem{prop}[thm]{Proposition}
\newtheorem{cor}[thm]{Corollary}
\newtheorem{proper}[thm]{Property}
\newtheorem{claim}[thm]{Claim}
\newtheorem{noname}[thm]{\ignorespaces}
\newtheorem*{thm1}{Theorem 1}
\newtheorem*{thm2}{Theorem 2}
\newtheorem*{thm1'}{Theorem $\mathbf{1'}$}
\newtheorem*{thm2'}{Theorem $\mathbf{2'}$}  
\theoremstyle{definition}
\newtheorem{defn}[thm]{Definition}    % Definition environment with 
\newtheorem{quest}[thm]{Question}
\newtheorem{notat}[thm]{Notation}
\newcommand{\hsp}{\hspace}  
\newcommand{\bitem}{\begin{itemize}}
\newcommand{\eitem}{\end{itemize}}
\newcommand{\ba}{\begin{array}}
\newcommand{\ea}{\end{array}}
\newcommand{\al}{\alpha}

\newcommand{\ga}{\gamma}

\newcommand{\om}{\omega}
\newcommand{\Om}{\Omega}

\newcommand{\lgl}{\langle}
\newcommand{\rgl}{\rangle}
\newcommand{\ra}{\rightarrow}

\newcommand{\lra}{\longrightarrow}

\newcommand{\hra}{\hookrightarrow}
\newcommand{\smetmi}{\smallsetminus}

\newcommand{\Int}{\int\limits}

\newcommand{\Sum}{\sum\limits}

\newcommand{\bcap}{\bigcap\limits}
\newcommand{\bcup}{\bigcup\limits}

\newcommand{\ol}{\overline}
\newcommand{\N}{\mathbb{N}} %{I\!\!N}

\newcommand{\R}{\mathbb{R}}

\newcommand{\cF}{{\cal F}}
\newcommand{\cG}{{\cal G}}

\begin{document}

\title{Foliations with few non-compact leaves}

\author{Elmar Vogt}                  
\address{2, Mathematisches Institut, Freie Universit\"at Berlin\\Arnimallee
3, 14195 Berlin, Germany}                  
\asciiaddress{2, Mathematisches Institut, Freie Universitaet Berlin\\Arnimallee
3, 14195 Berlin, Germany}

\email{vogt@math.fu-berlin.de}                     

\begin{abstract}  
Let $\cF$ be a foliation of codimension 2  on a compact manifold with at least
one non-compact leaf. We show that then $\cF$  must contain uncountably many
non-compact leaves. We prove  the same statement for
oriented 
$p$-dimensional foliations of arbitrary codimension if there
exists a closed  $p$ form which evaluates positively on every compact leaf.
For foliations of codimension 1 on compact manifolds  it is known that the 
union of all non-compact leaves is an open
set \cite{Hae}.

\end{abstract}
\asciiabstract{ 
Let F be a foliation of codimension 2 on a compact manifold with at
least one non-compact leaf. We show that then F must contain
uncountably many non-compact leaves. We prove the same statement for
oriented p-dimensional foliations of arbitrary codimension if there
exists a closed p form which evaluates positively on every compact
leaf.  For foliations of codimension 1 on compact manifolds it is
known that the union of all non-compact leaves is an open set [A
Haefliger, Varietes feuilletes, Ann. Scuola Norm. Sup. Pisa 16 (1962)
367-397].}

%\secondaryclass{  }              
\makeshorttitle  %%% Makes a short header for AGT and GTM articles

\setcounter{section}{-1}
\section{Introduction}
\label{intro}

Consider a $C^r$-foliation $\cF$ on a compact manifold with at least one
non-compact leaf. Is it possible that this leaf is the only non-compact leaf of
$\cF$? If not, is it possible that there are only finitely many non-compact
leaves, or countably many of them? Or must there always be uncountably many
non-compact leaves? Do the answers  depend on $r$? These
questions were asked by  Steve Hurder in
\cite{L},  Problem A.3.1.

 At first it seems obvious that for a foliation on a compact manifold
the union of all non-compact leaves, if not empty, should have a
non-empty interior. In fact, in codimension~1, and apart from flows,
these are the foliations that come first to mind, this set is open.
In \cite{Hae}, p.386, A. Haefliger proves that the union of all closed
leaves of a codimension~1 foliation of a manifold with finite first
mod~2 Betti number is a closed set.  Therefore, the union of all
non-compact leaves of a codimension~1 foliation of a compact manifold
is open. Consequently, the set of all non-compact leaves of a
codimension~1 foliation on a compact manifold is either empty or is
uncountable.

But for foliations of codimension greater than 1 it is easy to
construct examples on closed manifolds with the closure $C$ of the
union of all non-compact leaves a submanifold of positive
codimension. In fact, given $0<p<n$, there exist real analytic
$p$-dimensional foliations on closed $n$-manifolds where $C$ is a
$(p+1)$-dimensional submanifold.  Therefore, the dimension of the
closure of the union of all non-compact leaves can be quite small when
compared to the dimension of the manifold, even for real analytic
foliations (Compare this with Problem~A.3.2 in \cite{L} where a
related question is asked for $C^1$-foliations). These examples are
fairly straight forward generalizations of a construction of G. Reeb,
\cite{R},(A,III,c), and will be presented in Section~\ref{sec:dim}
(Proposition~\ref{prop:Reebrealanalytic}).

 The main results of this note extend the statement that for codimension~1
foliations on compact manifolds the set of non-compact leaves is either empty
or uncountable in two directions. First we show that it  also holds for
foliations of codimension~2. The second result states that it is true  in
general if an additional homological condition is satisfied. To be more
explicit, we first recall that a Seifert fibration is a foliation
whose leaves are all compact and  all have finite holonomy groups. Then we have:

\begin{thm1}
 Let $\cF$ be a foliation of codimension 2 on a compact
manifold. Then $\cF$ is either a Seifert fibration or it has uncountably
many non-compact leaves.
\end{thm1}

\begin{thm2} Let $\cF$ be an oriented $C^1$-foliation of dimension
$p$ on a compact manifold $M$. Assume that there exists a closed $p$-form
$\omega$ on
$M$ such that $\int_L \omega > 0$ for every compact leaf $L$ of $\cF$. Then $\cF$
is either a Seifert fibration or $\cF$ has uncountably many non-compact leaves.
\end{thm2}

Note that for a foliation on a manifold with boundary we will always assume
that the boundary is a union of leaves.

The two theorems are corollaries of farther reaching  but more
technical results stated further down in this introduction as theorems
$1^\prime$ and $2^\prime$.

We also include a short proof of the probably well known fact that
for an arbitrary suspension foliation (i. e. a foliated bundle) over a compact
manifold the set of non-compact leaves is empty or uncountable
(Proposition~\ref{prop:suspension}). 

As the statements of the two theorems indicate the techniques for their proofs
are strongly related to the methods in \cite{EMS} (and \cite{Vo1} for the
codimension-2 case). There it was shown that foliations with all leaves compact
on compact manifolds are Seifert fibrations if the homological
condition of theorem 2 holds or if the codimension is 2. The methods for the
codimension 2 case are essentially due to D.B.A. Epstein who proved the
corresponding result for circle foliations of compact 3-manifolds in \cite{Ep1}.
In these papers the result for foliations with all leaves compact and of
codimension~2 follows from the following more technical statement. Let $B_1$ be
the union of leaves with infinite holonomy. For a foliation with all leaves
compact of codimension 2 this set is empty if it is compact. This result in
turn is obtained by constructing a compact transverse 2-manifold $T$
intersecting each leaf of $B_1$ but with
$\partial T \cap B_1 = \emptyset$. One then uses a generalization of a theorem
of Weaver \cite{Wea} to show that there is a compact neighborhood $N$  of $B_1$
and an integer $n$ such that all but finitely many leaves of $N$ intersect $T$
in exactly $n$ points. Thus all holonomy groups of leaves in $N$ are finite and
$B_1 = \emptyset$. The construction of $T$ is by downward induction,
constructing transverse manifolds for a whole hierarchy $B_\alpha , \alpha$ an
ordinal, of so-called bad sets. Here, given $B_\al$, the set  $B_{\alpha+1}$
is defined as the union of all leaves of
$B_\alpha$ with infinite holonomy group when the foliation is restricted to
$B_\alpha$. (In \cite{Ep1}, \cite{EMS} and\cite{Vo1} a finer hierarchy of bad
sets is used. There a leaf $L$ of $B_\alpha$ belongs  to  $B_{\alpha+1}$
if the holonomy group of $L$ of the foliation restricted to $B_\alpha$ is
not trivial. The hierarchy we use in the present paper is called the coarse
Epstein hierarchy in \cite{EMS}.) 

For arbitrary codimension and still all leaves compact the ideas are due to
R.~Edwards, K.~Millett, and D.~Sullivan \cite{EMS}. They show that $B_1$ is
already empty if it is compact and if there exists  a closed form
$\omega$, defined in a neighborhood of $B_1$, satisfying $\int_L \omega > 0$ for
every leaf $L$ of $B_1$. The key idea in their proof is
to construct a sequence of homologous leaves in the complement of $B_1$
converging to $B_1$ such that the volume of the leaves grows to infinity as the
leaves approach $B_1$. Such a family gives rise to a non-trivial foliation cycle
which is roughly the limit of the leaves  each divided by some normalization
factor which tends to infinity as the volume of the leaves goes to
infinity. Therefore this cycle evaluates  on $\omega$ to 0 since 
integration of $\omega$ is constant on the homologous leaves. On the other hand,
if
$B_1$ is non-empty, the cycle cannot evaluate to 0 on $\omega$ since it is
non-trivial with support in $B_1$ and $\omega$ is positive when integrated over
any leaf of
$B_1$. (For a more detailed overview of this proof and an exposition of its
main ideas read the beautifully written introduction of \cite{EMS}).

In our situation we first extend
the notion of the hierarchy of bad sets  to incorporate the occurrence of
non-compact leaves. The points where the volume-of-leaf function with respect
to some Riemannian metric is locally unbounded, used in the papers mentioned
above for the definition of the first bad set $B_1$, is obviously inadequate.
Also the union of all leaves with infinite holonomy misses some irregularities
caused for example by simply connected non-compact leaves. Instead our
criterion in the inductive definition of the hierarchy of bad sets puts a leaf
of the bad set $B_\alpha$ in the next bad set $B_{\alpha+1}$ if for any
transversal through this leaf the number of intersection points with leaves of
$B_\alpha$ is not bounded. We begin with the whole manifold as $B_0$. In the
presence of non-compact leaves the first bad set
$B_{1}$ is not empty if the manifold is compact. 

As opposed to the case when
all leaves are compact, where the hierarchy of bad sets eventually reaches the
empty set, it is now possible that the hierarchy stabilizes at a non--empty bad
set
$B_{\al}$, i.~e.\
$B_{\al} = B_{\al+1}$ (and consequently $B_{\al} = B_{\beta}$ for all $\beta >
\al$) and $B_\al \neq \emptyset$.
But this will imply that $B_{\al}$ contains uncountably many
non--compact leaves (actually a bit more can be
said, see Proposition~\ref{prop:stabilization}). Thus, we may assume that the
hierarchy reaches the empty set. Then, in the codimension-2 case, we manage to
mimic all the steps in the construction of the transverse 2--manifold
$T$  mentioned above, if the following condition is satisfied: let
$N_{\al}$ be the union of the non-compact leaves in
$B_{\al} \smallsetminus B_{\al+1}$; then $\dim N_\al \leq\dim \cal{F}$. 
Thus, and by some (further) generalization of Weaver's theorem we obtain
the following theorem.

\begin{thm1'} Let $\cal{F}$ be a foliation of codimension 2,
let $B_{0} \supset B_{1} \supset \cdots$ be its Epstein hierarchy of
bad sets, and let $N_{\al}$ be the union of the non--compact leaves of
$B_{\al} \smallsetminus B_{\al + 1}$. If $B_{1} \cup N_{0}$ is compact,
then at least one of the following statements holds:
\begin{itemize}
\item[\rm(i)] for some ordinal $\al$ we have $B_{\al} = B_{\al+1}$ and
$B_{\al} \neq \varnothing$ -- in this case no leaf of $B_{\al}$ is
isolated (i.\ e.\ for no transverse open 2--manifold $T$ the set $T
\cap B_{\al}$ contains an isolated point) and $B_{\al}$ contains a dense
$G_{\delta}$ consisting of (necessarily uncountably many) non--compact leaves;
furthermore, $\dim B_\al > \dim \cF$, if all leaves of $B_\al$ are
non-compact \   --, or
\item[\rm(ii)] for some $\al~ \dim N_{\al} > \dim \cal{F}$, or
\item[\rm(iii)] $\cal{F}$ is a Seifert fibration.
\end{itemize}
\end{thm1'}

Theorem 1 follows from this since for any foliation  the set $B_{1} \cup N_{0}$
is closed.

That, in a way, Theorem 1$^\prime$ is best possible is shown by the examples of
Section~\ref{sec:dim} mentioned above. They contain examples of real
analytic foliations of codimension 2  on closed manifolds of any given
dimension greater than 2 such that $B_1$ consists of finitely many compact
leaves, and the dimension of the union of the non-compact leaves exceeds the
leaf dimension by one.

The procedure for the proof of Theorem
2 is also similar to the proof in the case where all leaves are assumed to be
compact. The sequence of homologous compact leaves for the definition of the
foliation cycle is now assumed to converge to  $B_1 \cup N_0$ where as
above $N_0$ is the union of all non-compact leaves in the complement of $B_1$,
and the closed form $\omega$ has to be defined in a neighborhood of the closed
set
$B_1 \cup N_0$. The construction of this sequence of leaves is to a certain
degree easier in the presence of (not too many) non-compact leaves. Also the
support of the limiting foliation cycle will essentially be disjoint from the
non-compact leaves and thus they play no role in the evaluation of this cycle
on 
$\omega$. More precisely, we have the following theorem.

\begin{thm2'} Let $\cF$ be an oriented $C^1$-foliation of
dimension $p$ on a manifold $M$, let $M = B_{0} \supset B_{1} \supset \cdots$
be its Epstein hierarchy of bad sets, and let $N_0$ be the union of the
non--compact leaves in the complement of $B_{1}$. Assume that there
exists a closed $p$-form $\omega$ defined in a neighborhood of $B_1 \cup N_0$
such that $\int_L \omega > 0$ for every compact leaf $L$ of $B_1$, and assume
that $B_1 \cup N_0$ is compact. Then at least one of the following
statements holds:
\begin{itemize}
\item[\rm(i)] $ \bigcap_\alpha B_\alpha \neq \emptyset$ -- in this
case  no leaf of $\bigcap_\alpha B_\alpha$ is isolated and $\bigcap_\alpha
B_\alpha$ contains a dense $G_\delta$--set consisting of (necessarily
uncountably many) non--compact leaves --, or  
\item[\rm(ii)] $N_0$ is a non-empty open subset of $M$,
in fact a non-empty union of components of the open set $M\setminus B_1$, or   
\item[\rm(iii)] $\cF$ is a Seifert fibration.  
\end{itemize}
\end{thm2'}

 In all the examples that I am aware of where a $p$-dimensional foliation
on a compact manifold contains non-compact leaves, the union of all non-compact
leaves is at least $p+1$ dimensional.  Therefore, it would be interesting to
know whether statement {(i)} in the two theorems above could be improved
to: 
$B_1$ contains a subset of dimension greater than the leaf dimension consisting
of non-compact leaves. I also do not know whether in codimension $\geq 3$ there
are foliations on compact manifolds with at least 1 and at most countably many
non-compact leaves.  

Theorems 1 and 1$^\prime$ hold for topological foliations, but we give a
detailed proof only for the $C^1$--case, indicating the necessary changes for the
topological case briefly at the end of section~\ref{sec:codim2}. The proof in
the
$C^0$--case depends heavily on the intricate results of D.B.A. Epstein in
\cite{Ep3}.

Section~\ref{sec:dim} contains the examples mentioned above of
real analytic foliations of codimension $q$, $q\geq1$, on closed manifolds such
that the closure of the union of all  non--compact leaves is a submanifold of
codimension $q-1$. 

If one is content with $C^{r}$--foliations, $0\leq r\leq\infty$, then
one can construct such examples on many manifolds: given $p, k,n$ with $0<p < k
\leq n$ and $p\leq n-2$ then any
$n$-manifold which admits a $p$-dimensional $C^{r}$-foliation with all leaves
compact will support also a $p$-dimensional $C^r$-foliation such that the
closure of the union of all non-compact leaves is a non-empty submanifold of
dimension
$k$\   (Proposition~\ref{prop:Reeb}).

In addition, we give in Section~\ref{sec:dim} a simple proof of the well-known
fact that for suspension foliations over compact manifolds, i. e. foliated
bundles with compact base manifolds, the existence of one non-compact leaf
implies the existence of uncountably many. This is an easy application of a
generalization, due to D.B.A. Epstein \cite{Ep2}, of a theorem of Montgomery
\cite{M}.

In Section~\ref{sec:isolated} we introduce some notation and gather a few
results concerning the set of non--compact leaves of a foliation. In
particular, we prove a mild generalization of the well known fact that 
the closure of a non-proper leaf of a foliation contains uncountably many
non-compact leaves.

 Section~\ref{sec:hierarchy} introduces the notion of
Epstein hierarchy of bad sets in the presence of non--compact leaves and we
prove some of its properties. Section~\ref{sec:codim2} contains the proof of
Theorem 1$^\prime$ along the lines indicated above. Finally, in
Section~\ref{sec:codimp} we construct the sequence of compact leaves
approaching
$B_1 \cup N_0$ and the associated limiting foliation cycle, and establish the
properties of this cycle to obtain the proof of Theorem~2$^\prime$.  

I have tried to make this paper reasonably self contained. But
in Section~\ref{sec:codim2} referring to some passages in \cite{Vo2} will
be necessary for understanding the proofs in all details. The same holds
for Section~\ref{sec:codimp} where familiarity with \cite{EMS} will be very
helpful. I will give precise references wherever they are needed. Also, I will
use freely some of the notions and results of the basic paper \cite{Ep2} on
foliations with all leaves compact.

It will also be of help to visualize some of the examples of
Section~\ref{sec:dim}. Although they are simple they illustrate some of the
concepts introduced later in the paper, and they give an indication of the
possibilities expressed in Theorems $1'$ and $2'$ above.

This paper replaces an
earlier preprint of the author with the same title. There the main result of
\cite{Vo3} and calculations of the Alexander cohomology of the closure of
the union of all non-compact leaves was used to prove a special case of Theorem
1 for certain 1-dimensional foliations on compact 3-manifolds. 

In this paper finite numbers are also considered to be countable.

\section{Foliations having a set of
non--compact leaves of  small dimension}
\label{sec:dim}

We generalize (in a trivial way) an example given by G.~Reeb in
\cite{R},(A,III,c).  Let $F^{p}$ and $T^{k}$ be closed connected real analytic
manifolds. Let $f: F^{p}\ra\R$ be real analytic with $0$ a regular value in
the range of $f$, and let
$g: T^{k}\ra\R$ be real analytic with a unique maximum in $x_0 \in T^k$.  For
convenience, let
$g(x_{0}) =1$.

Let $\theta\in\R\,\mathrm{mod}\,2\pi$ be coordinates for $S^{1}$ and consider
for
$x\in T^{k}$ the 1--form
$$
\omega(x) = ((g(x)-1)^{2} +f^{2}) d\theta + g(x) df
$$
on $F^{p}\times S^{1}$.  One immediately checks that $\omega(x)$ is
nowhere $0$ and completely integrable. It thus defines a real analytic foliation
$\cF(x)$ of codimension~1 on $\{x\} \times F^{p} \times S^{1}$. 
These foliations fit real analytically together to form a foliation
$\cF$ of codimension $q=k+1$ on $T^{k} \times F^{p} \times S^{1}$. 
It is easy to describe the leaves of $\cF(x)$.  There are two cases.

{\bf Case 1}\qua  $x\neq x_{0}$.  Then $(g(x)-1)^{2} >0$, and 
$\omega(x) =0$ if and only if $\displaystyle d\theta =
-\frac{g(x)}{(g(x)-1)^{2}+f^{2}}\, df$. Therefore the leaves of $\cF(x)$ are the
graphs of the functions $h(\theta_{0}): F^{p}\ra S^{1}$, given by
$$ h(\theta_{0}) (y) = \frac{-g(x)}{g(x)-1}\cdot \arctan
\big(\frac{f(y)}{g(x)-1}\big) + \theta_{0}, \, 0\leq \theta_{0} < 2\pi.$$  
These leaves are all diffeomorphic to $F^{p}$ and therefore compact.

{\bf Case 2}\qua  $x=x_{0}$.  Then $g(x)=1$, and $\omega(x_{0}) = f^{2}
d\theta +df$.  We obtain two kinds of leaves for $\cF(x_{0})$.  Let
$F_{0} = f^{-1}(0)$. Then $\omega(x_{0}) = df$ on $\{x_0\}\times F_{0} \times
S^{1}$.  Therefore, $\omega(x_{0}) =0$ implies $f=\mbox{ const}$, and
the components of $\{x_0\}\times F_{0} \times S^{1}$ are compact leaves of
$\cF(x_{0})$.  In the complement of $\{x_0\}\times F_{0}\times S^1$ the
foliation
$\cF(x_{0})$ is given by $\displaystyle d\theta = -\frac{1}{f^{2}}\, df$. 
Therefore, the leaves are components of the graphs of $k(\theta_{0}): F\smetmi
F_{0} \ra S^{1}$, given by
$$
k(\theta_{0}) (y) =  \frac{1}{f(y)} +\theta_{0}, \, 0\leq\theta_0 < 2\pi.
$$
Since $F$ is connected, no component of $F-F_{0}$ is compact.  Thus
$\{x_{0}\} \times (F\smetmi F_{0}) \times S^{1}$ is the union of the 
non--compact leaves of $\cF$, and we have the following result.

\begin{prop}
\label{prop:Reebrealanalytic} Let $F^{p}$ and
$T^{k}$ be real analytic closed manifolds of dimension $p > 0$ and $k$
respectively.  Then there exists a real analytic $p$-dimensional foliation
$\cF$ on $T^{k} \times F^{p} \times S^{1}$ such that the closure of the union
of the non--compact leaves of $\cF$ equals $\{x_{0}\} \times F^{p}
\times S^{1}$ for some point $x_{0}\in T^{k}$. \qed
\end{prop}

The construction above is quite flexible, especially if one allows
the foliations to be smooth.  For example

\begin{prop}
\label{prop:Reeb} Let $0\leq r\leq \infty$ and let $M$ be an $n$-manifold
which supports a $p$-dimensional
$C^{r}$--foliation with all leaves compact, where $0<p\leq n-2.$  Then for
any integer $k$ with $p<k\leq n $ the manifold $M$ supports a p-dimensional
$C^{r}$--foliation with the following property:  the closure of the union of
non--compact leaves is a  non-empty submanifold of dimension $k$.
\end{prop}

\begin{proof}  The leaves with trivial holonomy form an open dense
subset of any foliation with all leaves compact.  Let $F$ be a leaf
with trivial holonomy and $U$ a saturated neighborhood of $F$ of the
form $F\times D^{n-p}$, foliated by $F\times \{y\}$, $y\in D^{n-p}$,
where $D^{n-p}$ is the unit $(n-p)$--ball.  Let $S^{1}\times D^{n-p-1} \hra
D^{n-p}$ be a smooth embedding into the interior of $D^{n-p}$ and let $K$ be a
compact submanifold of the interior of $D^{n-p-1}$ of dimension $k-p-1$. 

Let $f: F\ra \R$ be smooth with $0$ a regular value in the range of $f$, 
and let $h: D^{n-p-1}\ra [0,1]$ be  smooth  with
the following properties:\\
$h$ and all its derivatives vanish on
$\partial D^{n-p-1}$, $h(K)=1$, and $h(z)<1$ for all $z\in
D^{n-p-1}\smetmi K$.

Replace on
$F\times S^{1}\times D^{n-p-1}$ the product foliation induced from $F\times
D^{n-p}$ by the smooth foliation defined on $F\times S^{1}\times \{z\}$, $z\in
D^{n-p-1}$, by the 1--form
$$
\omega(z) = ((h(z)-1)^{2} + f^{2}) d\theta + h(z) df \,.
$$
Then, as in our first example, the foliations of codimension~1 on $F\times
S^{1}\times\{z\}$ defined by
$\omega(z) =0$ fit smoothly together to form a foliation of $F\times
S^{1}\times D^{n-p-1}$. On the boundary $F\times S^{1}\times
\partial D^{n-p-1}$  this foliation fits smoothly to the product foliation on
$F\times (D^{n-p}
\smetmi (S^{1}\times D^{n-p-1}))$.  The leaves in $F\times S^{1}\times
(D^{n-p-1}\smetmi K)$ are all diffeomorphic to $F$.  Furthermore,
$(F\smetmi f^{-1}(0)) \times S^{1} \times K$ is a union of
non--compact leaves, each one of which is diffeomorphic to some
component of $F\smetmi f^{-1}(0)$.  Therefore $F\times S^{1}\times
K$ is the closure of the union of the non--compact leaves. 
\end{proof}

Proposition~\ref{prop:Reeb} suggests the following question.

\begin{quest}
\label{quest:dim} Does there exist a foliation of dimension 
$p$ on a compact manifold $M$ such that the closure of the union of all
non--compact leaves is non--empty and has dimension $p$\,?
\end{quest}

By the result of Haefliger mentioned in the introduction \cite{Hae}, page 386,
such a foliation has codimension at least 2, and in the case of codimension
2, there must, by Theorem $1'$, be an $\al$ such that $ B_\al = B_{\al+1}\neq
\emptyset$, and $B_\al$ contains compact leaves. 
Furthermore, in general, it cannot be a suspension foliation.  A suspension
foliation
$\cF_{\varphi}$ is given by a homomorphism 
$\varphi:
\pi_{1}(B) \ra \mathrm{Homeo}(T)$, where $B$ and $T$ are 
manifolds and $B$ is connected.  One foliates $\tilde{B}\times T$ by
$\tilde{B}\times
\{t\}$, $t\in T$, where $\tilde{B}$ is the universal cover of $B$. 
This foliation is invariant under the obvious action of $\pi_{1}(B)$
and induces the foliation $\cF_{\varphi}$ on the quotient $\tilde{B}
\times T/\varphi$ by this action.

\begin{prop}
\label{prop:suspension}
 Let $B$ and $T$ be  manifolds with 
$B$ closed and connected, let  \mbox{$\varphi: \pi_{1}B \ra  Homeo(T)$} be a
homomorphism and assume that the associated suspension foliation
$\cF_{\varphi}$ has non--compact leaves. Let $N$ be the closure of the union
of the non-compact leaves, and let $W$ be a component of $(\tilde{B}
\times T/\varphi) \; \setminus N$. Then the closure of $W$ consists of
compact leaves.  In particular, $N$ does not contain any isolated
leaf, and $(\tilde{B} \times T/\varphi) \; \setminus N$ consists of infinitely
many components unless N contains interior points. Furthermore, the
dimension of $N$ is at least equal to $( \dim B+\dim T-1)$. (As definition
of dimension we may take any of the notions of covering dimension, inductive
dimension, or cohomological dimension, which are equivalent in our situation).

If $\dim T=1$, then the union of all non-compact leaves is open.
\end{prop}

\begin{proof}  The quotient space $M = \tilde{B}\times T/\varphi$ is a
fibre bundle with fibre $T$.  The fibres are transverse to the
foliation $\cF_{\varphi}$.  Compact leaves of $\cF_{\varphi}$
correspond to finite orbits of the group $G = \varphi(\pi_{1}(B)) \subset
\mbox{Homeo}(T)$. We identify $T$ with the fibre over the basepoint of $B$.  
Let $W$ be a component of $M \setminus N$ and let $x$ be a point of $W \cap T$.
Let $G_{Ox}$ be the normal subgroup of finite index of $G$ whose elements keep
the orbit of $x$ pointwise fixed. Let $W_x$ be the component of $W \cap T$
containing $x$ and let $G_W$ be be the restriction of $G_{Ox}$ to $W_x$
considered as a subgroup of
$\mbox{Homeo}(W_x)$. Every orbit of $G_W$ is finite and $W_x$ is a connected
manifold. Then, by Theorem~7.3 in \cite{Ep2}, an extension of the main result
in \cite{M} to groups of homeomorphisms, 
$G_W$ is finite, say of order $g$. This implies that the $G_{Ox}$-orbit of any
point in the closure of $W_x$ contains at most $g$ points. Thus any orbit
in the closure of
$W
\cap T$ under the action of
$G$ is finite, and  all leaves in the closure of
$W$ must be compact. 

Since this is true for any component $W$ of $M \setminus N$ any
neighborhood of a point in
$N$ which is not a point of the closure of $\mathrm{int}(N)$ must intersect
infinitely many components of
$M\setminus N$. In particular, $\dim N \geq \dim M-1= \dim B+\dim T-1$. 

If $\dim
T=1$, we first reduce our problem to the case where $T$ is connected. To do this we
observe that $T$ decomposes naturally into disjoint subspaces each of which
is a union of components of $T$ on which $G$ acts transitively.  So we may assume
that $G$ acts transitively on the components of $T$.  If then the number of
components is infinite, then all leaves are non-compact. So we may assume that the
number of components of
$T$ is finite. We then replace $G$ by the subgroup of finite index whose elements
preserve every  component. This corresponds to passing to a finite covering of $M$
with the induced foliation. On each component of this covering space the induced
foliation is a suspension with one component of $T$ as fibre. 
So we may assume that $T$ is connected, i.e. $T$ is either $\mathbb{R}$ or $S^1$. 
We may furthermore assume that all elements of
$G$ are orientation preserving.  

If $T=\mathbb{R}$, then every finite orbit is
a global fixed point. Therefore, the union of all finite orbits is closed, and
we are done. If $T=S^1$ either all leaves are non-compact or we may pass to
the subgroup of finite index keeping a finite orbit pointwise fixed. For this
subgroup every finite orbit is again a global fixed point, and we can argue as
before. 
\end{proof}

\begin{cor}[Well known]
\label{cor:suspension} Let $\cF_{\varphi}$ be the
suspension foliation associated to the homomorphism $\varphi:
\pi_{1}(B)\ra \mathrm{Homeo}(T)$ and assume that $B$ is closed.  Then
$\cF_{\varphi}$ contains uncountably many non--compact leaves or none. If $\dim
T = 1$, the union of all non-compact leaves is open.
\end{cor}

\begin{proof}  Assume that $\cF_{\varphi}$ contains a non--compact
leaf, and let $N$ be the closure of the union of the non--compact
leaves of $\cF_{\varphi}$. The set
$N$ with its induced foliation is a foliated space in the sense of \cite{EMT}. 
The main result of \cite{EMT} implies that the union of all leaves of $N$ with
trivial holonomy is a dense
$G_{\delta}$ in $N$.  By Proposition~\ref{prop:suspension},  $\dim N\geq
\dim\cF_{\varphi} +1$ and $N$ does not contain any isolated leaf. Therefore, the
Baire category theorem implies that $N$ contains uncountably many leaves with
trivial holonomy.  Let $L$ be a compact leaf of $N$.  Then every neighborhood
of $L$ intersects a non--compact leaf (of $N$).  Therefore, the holonomy of $L$
is non--trivial, i.e. the leaves of $N$ with trivial holonomy are all
non--compact.
\end{proof}

\section{Uncountably many versus
isolated non--compact leaves}
\label{sec:isolated}

The material in this short section is  standard.  We include it to
fix and introduce notation.

Let $\cF$ be a foliation of a manifold $M$ and let $A\subset M$ be a
union of leaves.  We call a leaf $L\subset A$ {\it isolated (with
respect to $A$)} if $L$ is an open subset of $Cl(A)$.  Recall that a
leaf of a foliation is called proper if its leaf topology coincides
with the induced topology as a subset of $M$.  Obviously, any leaf
which is isolated with respect to some $A$ is proper, and the proper
leaves are exactly those leaves which are isolated with respect to
themselves.  We will denote the union of isolated leaves with respect
to the saturated set $A$ by $I(A)$.

\begin{prop} Let $A$ be a union of non--compact 
leaves of a foliation.  Then at least one of the following holds:
\bitem
\item[\rm (i)]  $Cl(A)$ contains uncountably many non--compact leaves,
or 
\item[\rm (ii)]  $A\subset Cl(I(A))$.
\eitem
In particular, the closure of a non-proper leaf contains uncountably many
non-compact leaves.
\end{prop}

\begin{proof}  Assume that $B:= A\smetmi Cl(I(A))$ is not empty.  Any
isolated leaf with respect to $B$ is also isolated with respect to
$A$.  Therefore $I(B) =\emptyset$.  Let $U = D'\times D$ be a
foliation chart with $D'$ connected and tangent to the foliation and
assume that $U\cap B\neq \emptyset$. We identify $D$ with $\{y\}\times D$
for some basepoint $y\in D'$.   Then
$C := Cl(B)
\cap  D$ is a closed non--empty subset of $D$ which
contains no isolated points.  By the main theorem of \cite{EMT} the union
$H$ of leaves of $Cl(B)$ with trivial holonomy is a dense
$G_{\delta}$. Since all leaves in $B$ are non-compact, the compact leaves in
$Cl(B)$  all have non-trivial holonomy. Therefore all leaves in
$H$ are non--compact.  Since
$H\cap C$ is a dense $G_{\delta}$ in $C$, the set $H\cap C$ is uncountable by
the Baire category theorem.  Every leaf of our foliation intersects
$D$ in an at most countable set.  Therefore the set of
non--compact leaves of $Cl(B)$ intersecting $D$ is uncountable. 
\end{proof}

\section {The Epstein hierarchy in
the presence of non--compact leaves}
\label{sec:hierarchy}

There are several possible ways  to generalize the notion of Epstein
hierarchy (see \cite{Ep1} or \cite{Vo2}) to foliations admitting
non--compact leaves. For our purposes definition~\ref{defn:hierarchy} below
seems to be the best choice. Before we come to this we need some notation.

\begin{notat}
\label{not:transverse} Let $(M, \cF)$ be a codimension $k$ foliated
manifold. A transverse manifold is a $k$-dimensional submanifold $T$
of $M$ which is transverse to $\cF$ and whose closure $Cl(T)$ is
contained in the interior of a $k$-dimensional submanifold transverse
to $\cF$. A transverse manifold $T$ may or may not have a boundary, 
denoted by $\partial T$. We call $\mathrm{int} T : = T \smetmi \partial T$ 
the interior of $T$, and call $T$ open, if $T=\mathrm{int} T$.
\end{notat}

\begin{notat}
\label{not:sec} Let $T$ be a transverse manifold of a
foliated manifold $(M, \cF)$. Then $\sec_{T} : M \lra \ol{\N} : = \N
\cup \{\infty \}$ is the map which associates to $x \in M$ the
cardinal of the set $T \cap L_{x}$, where $L_{x}$ is the leaf through
$x$.
\end{notat}

The topology on $\ol{\N}$ is the one point compactification of $\N = \{ 0,
1, 2, \ldots \}$. Then we have

\begin{proper}
\label{proper:sec}
\bitem
\item[\rm (a)] {\it $\sec_{T}$ is continuous in every point of $\sec^{-1}_{(T
\smetmi \partial T)} (\infty)$.}
\item[\rm (b)] {\it If $L_{x} \cap \partial T = \emptyset$, then
$\sec_{T}$ is lower semi-continuous in $x$.}
\eitem 
\end{proper}

The proofs are obvious. \qed

\begin{defn}
\label{defn:hierarchy} Let $(M, \cF)$ be a foliated manifold. The
Epstein hierarchy of bad sets of $\cF$ is a familiy $\{B_{\al} =
B_{\al} (\cF) \}$ of subsets of $M$ indexed by the ordinals and is
defined by transfinite induction as follows:
$$
\begin{array}{lclll}
B_{0} & = & M & ;\\
B_{\al} & = & \bcap_{\beta < \al} B_{\beta} & , & \textrm{if }  \al
\textrm{ is a limit ordinal;}\\
B_{\al + 1} & = & \{x \in B_{\al} & : & \textrm{for every transverse
manifold\  } T \textrm{\  with } x \in \mathrm{int}T\\
& & & & \sup \{\sec_{T} (y) : y \in B_{\al}\} = \infty\}~.
\end{array}
$$
\end{defn}

\noindent
Obviously each $B_{\al}$ is a closed invariant set.

\begin{prop}
\label{prop:stabilization} \bitem
\item[\rm(1)] If $B_{\al + 1} = B_\al$ and $B_\al \neq
\emptyset$ then for any transverse open manifold $T$ with $T\cap B_\al \neq
\emptyset$ we have

\bitem
\item[\rm (i)] $T\cap B_\al$ contains no isolated point and
\item[\rm (ii)] $T\cap B_\al$ contains a dense (necessarily uncountable) 
$G_\delta$-set
$R$ of points lying in non-compact leaves.
\eitem 

\item[\rm (2)]If $B_{\al + 1}$ contains at most countably many
non-compact leaves, then $B_{\al + 1}$ is nowhere dense in $B_{\al}$.
\eitem
\end{prop}

\begin{proof} By the Baire
category theorem a locally compact space without isolated points does not
contain a countable dense $G_\delta$-set. For any transverse open manifold $T$
the set $T\cap B_\al$ is locally
compact, and its isolated points belong to $B_\al \setminus B_{\al + 1}$. Any
leaf intersects any transverse manifold in an at most countable set. Therefore
it suffices to show the following. For every open transverse manifold
$T$ with $T\cap B_\al \neq \emptyset$ and $T\cap B_\al \subset B_{\al +1}$ the
set $T\cap B_\al$ satisfies properties (i) and (ii). 

Let sec$_\al$ be the restriction of sec$_T$ to the locally compact space
$T\cap B_\al$. By (\ref{proper:sec}) and again the Baire category theorem $T
\cap B_{\al}$  contains a dense $G_{\delta}$-set $R$ of points where
$\sec_{\al}$ is continuous. Assume that there exists a point  $y
\in R \cap \sec^{-1}_{\al} (\N)$. Then $\sec_{\al}$ is  constant in
a neighborhood of  $y$. This means that $y$ is a point in $B_{\al} \smetmi
B_{\al + 1}$, which is not possible. Therefore, $R \subset
\sec^ {-1}_{\al} \{\infty \}$. This implies that every point of $R$ is contained
in a non-compact leaf of $B_{\al + 1}$. 
\end{proof}

The  first claim of the next proposition is due to the convention that manifolds
are second countable. 

\begin{prop}
\label{prop:hierarchyend}
If $B_1\neq B_0$, let $\gamma:= \min\{\beta \;|\;
B_{\beta +1}= B_{\beta +2}\}$. Otherwise, let $\gamma$ be $0$. Then the
following holds:
\bitem
\item[\rm (i)]   $\gamma$ is a countable ordinal;
\item[\rm (ii)]  if $M\neq \emptyset$ and $B_{\al}$ is compact for some $\al <
\ga$ then $B_{\ga}
\neq
\emptyset$;
\item[\rm (iii)] if
$B_{\ga}$ is compact then all leaves in $B_{\ga}$ are compact.
\eitem
\end{prop}

The last statement is due to the fact that $B_{\beta + 1} \neq \emptyset$ if
$B_\beta$ is compact and contains a non-compact leaf.
\bigskip

For further reference we note the following proposition.

\begin{prop}
\label{prop:openbadset} Each point of the interior of $B_1$ is contained in
$\displaystyle \bigcap_\al B_\al$.
\end{prop}

\begin{proof} Let $x$ be an interior point of $B_1$, and let $T$ be an open
transverse manifold with $x\in T\subset B_1$. Then sec$_T$ = sec$_{T\cap B_1}$.
Thus $T\subset B_2$, and, by transfinite induction, $T\subset B_\al$ for all
$\al$. 
\end{proof}

\section{Codimension 2 foliations}
\label{sec:codim2}

For simplicity we assume that all foliations are $C^{1}$, but the
main result (Theorem 1$^\prime$) is also true for $C^{0}$-foliations. We
will indicate the necessary changes in an appendix at the end of this section.

As in the case of the study of foliations with all leaves compact
there are two ingredients which make the codimension 2 case special.
The first one is the fact that for $\al \geq 1$ the bad set $B_{\al}$ is
transversally of dimension at most 1 if $\dim\bigcap_\al B_\al < \dim B_0$
(Proposition~\ref{prop:openbadset}). The second one is a generalization of
Weaver's Lemma \cite{Wea} which takes in our setting the following form.

\begin{prop}
\label{prop:Weaver} Let $\cF$ be a foliation of codimension
2 and $T$ a transverse 2-manifold. Let $C \subset T$ be compact
connected and $W$  be the union of all leaves through points of $C$. Let $E
\subset C$ be the set of points of $C$ lying in a non-compact leaf. We assume
that

\bitem
\item[{\rm (i)}]  no compact leaf of $W$ intersects $\partial T$, 
\item[{\rm (ii)}]  every non-compact leaf of $W$ intersects $T$ in
infinitely many points,
\item[{\rm (iii)}]  for any loop $\om$ of a compact leaf through a point $x \in
C$ a representative of the associated holonomy map defined in a
neighborhood of $x$ in $T$ preserves the local orientation, and
\item[{\rm (iv)}]  $E$ is a  countable union of disjoint closed sets $E_j$.
\eitem

\noindent
Then either all leaves of $W$ are non-compact, or there exists an integer $\rho$
such that all but
finitely many leaves of $W$ intersect $T$ in
exactly $\rho$ points. In the latter case the finitely many other leaves of $W$
intersect $T$ in fewer than $\rho$ points.
\end{prop}

\begin{proof} For each positive integer $m$ let
$$C_{m} = \{x \in C : \sec_{T} (x) \leq m \},$$
and let $D_{m} \subset C_{m}$ be the set of non-isolated points of
$C_{m}$. Clearly each $C_{m}$ and $D_{m}$ is closed, each $C_{m}
\smetmi D_{m}$ is at most countable, and we have a
decomposition
$$C = \bcup_{m \geq 1} (D_{m} \smetmi D_{m-1}) ~~ \cup ~~
\{\mbox{countable set}\} \cup \bcup_{j} E_j$$
\noindent
into a countable union of disjoint sets. We claim that each $D_{m}
\smetmi D_{m-1}$ is closed. For if not, then there exists $x \in
D_{m-1}$ with $x \in Cl(D_{m} \smetmi D_{m-1})$. By hypothesis (ii)
the leaf $L_{x}$ through $x$ is compact. Now we can argue as in the
proof of Lemma 3.4 in \cite{Vo2} to obtain a representative h of an element
of the holonomy group of $L_{x}$ such that $dh(x)$ has a non-zero fixed
vector $v$ and a periodic vector $w$ with least period $\nu > 1$. But this
contradicts hypothesis (iii). (Hypothesis (i) is needed for imitating the proof
of 3.4 in \cite{Vo2}.)

Since by hypothesis (iv) also the $E_j$ are closed, the compact connected set
$C$ is a countable disjoint union of closed sets. Then a theorem of Sierpinski
\cite{Ku}, \S 47 III Theorem 6, states that $C$ must be equal to one of
the sets of which it is the disjoint union. So either $C$ is one of the $E_j$
and all leaves of $W$ are non-compact, or $C$ is a single point in a compact
leaf, or  there exists $\rho$ with $C = D_{\rho} \smetmi
D_{\rho-1}$.  The set of points of $ D_{\rho} \smetmi
D_{\rho-1}$ which lie in leaves intersecting $T$ in less than $\rho$ points
is $( D_{\rho} \smetmi D_{\rho-1})\cap C_{\rho - 1}$. But this set is compact
and discrete and therefore finite. 
\end{proof}

An easy consequence of 4.1 is the following result.

\begin{prop}
\label{prop:notransversal} Let $\cF$ be a codimension 2 foliation,  $B_{0}
\supset B_{1} \supset \cdots$ the Epstein hierarchy of $\cF$, and $N_\al$ the
union of all non-compact leaves of $\cF$ in $B_\al\setminus B_{\al+1}$. Assume
that the closure of every leaf of $\cF$ is compact, that $B_1$ is a non-empty
set, and that  $\dim N_\alpha\leq\dim\cF$ for every $\alpha\geq0$. Furthermore
assume that $\bigcap_\al B_\al$  
is empty or  contains only non-compact leaves  and that $ \dim
\bigcap_\al B_\al <\dim B_0$. Then there does not exist a compact
transverse 2-manifold
$T$ intersecting each leaf of
$B_{1}$ and with
$\partial T \cap B_{1} = \emptyset$.
\end{prop}

\begin{proof}  The proof is by contradiction. It is clear that we may assume
that
$\cF$ is transversely orientable. We will show below
(Lemma~\ref{lem:avoidboundary}) that with our hypotheses we can always arrange
$T$ so that $\partial T$ does not intersect any non-compact leaf. The
union of non-compact leaves in $B_{0} \smetmi B_{1}$ is closed in 
$B_{0} \smetmi B_{1}$. Therefore we find a compact neighborhood $K$ 
of $\partial T$ in $T$ such that $K \cap B_{1} = \emptyset$ and
every leaf through a point of $K$ is compact. This implies that the
union $S_{K}$ of leaves through $K$ is compact and $\cF$ restricted
to $S_{K}$ is a Seifert fibration (Here we extend the notion of Seifert
fibration to foliated sets. Such a set will be called a Seifert fibration if all
leaves are compact with finite holonomy groups).

Now consider a component $D$ of $T \smetmi S_{K}$ such that $D \cap
B_{1} \neq \emptyset$ (such a component exists) and apply
Proposition~\ref{prop:Weaver} to $C = \ol{D} \subset T \setminus {\rm Int}(K)$.
Hypotheses (i), (ii) and (iii) of \ref{prop:Weaver} are clearly satisfied, the
last one because we have assumed $\cF$ to be transversely orientable.
The union $N_\al$ of all non-compact leaves in
$B_\al\setminus B_{\al+1}$ is a closed subset of $B_0 \setminus B_{\al+1}$
which intersects $T$ in a set of dimension  $0$. Therefore, for each $\al$, the
set $T\cap N_\al$ is a countable disjoint union of compact sets.
By Proposition~\ref{prop:hierarchyend} there are only countably many
non-empty $N_\al$'s.  Furthermore, if $\bcap B_\al$ contains non-compact
leaves, all leaves of $\bcap B_\al$ are non-compact. Since $\bcap B_\al$  is
closed, also Condition (iv) is satisfied, and  we are
entitled to apply
\ref{prop:Weaver}. Since $D$ is a non-empty open subset of
$T$ and, by hypothesis, the union of all non-compact leaves has dimension less
than the dimension of $B_0$, not all points of
$\ol{D}$ lie in non-compact leaves. Consequently all leaves intersecting
$\ol{D}$ are compact and the function
$\sec_{T}$ is bounded on $\ol{D}$. This implies $B_{1} \cap {D} = \emptyset$,
which is a contradiction. 
\end{proof}

The next (easy) lemma is true in by far more generality. We only state it for
the case of interest to us.

\begin{lem}
\label{lem:avoidboundary} Let $T$ be a 2-manifold with compact boundary 
$\partial T$ and let  $N \subset T$ be a 0-dimensional subset which is closed
in a neighborhood of $\partial T$. Then for any neighborhood $U$ of $\partial T$
there exists a submanifold $T^\prime \subset T$ with compact boundary $\partial
T^\prime$ such that $T \setminus U \subset T^\prime $ and $\partial T^\prime
\cap N = \emptyset.$
\end{lem}

\begin{proof} By looking at each component of $\partial T$ separately the lemma
reduces to the statement that for any closed 0-dimensional subset $N$ of $S^1
\times [0,1]$ we find a neighborhood $K$ of $S^1 \times \{0\}$ which is a
compact 2-manifold with boundary $\partial K$ such that $S^1 \times
\{0\}\subset \partial K$ and 
$\partial K \cap N = (S^1 \times \{0\}) \cap N$.

Since $N$ is 0-dimensional and closed, $N$ is for any $\epsilon > 0$ a finite
disjoint union of closed sets of diameter less than $\epsilon $, where  we
metrize  $S^1 \times [0,1]$  by considering it as a smooth submanifold of $R^2$.
In particular, $N$ is the union of closed sets $N_0 $ and $N_1$ with
$(S^1 \times \{0\} \cup N_0) \cap  (S^1 \times \{1\} \cup N_1) = \emptyset$.  Let
$d$ be the distance between $(S^1 \times \{0\}) \cup N_0$ and  $(S^1 \times
\{1\})
\cup N_1.$  Then there exist finitely many closed disks $D_1, D_2, \ldots ,
D_s$ of radius $r < d$ such that $\{{\rm Int}D_i\}$ covers $(S^1 \times \{0\}
\cup N_0)\}$, their boundaries $\{\partial D_i\}$ are in general position, and
$D_i
\cap N_1 = \emptyset$ for all $i$. Then $K := \bigcup _i D_i\cap (S^1 \times
[0,1])$ is a 2-manifold with  piecewise smooth boundary having the desired
properties. We may, if we want to, smooth $\partial K$. Taking the component
containing $S^1
\times \{0\}$ and filling in some components of $\partial K$ bounding 2-cells
in $S^1 \times [0,1]$ we may also assume that $K$ is an annulus. 
\end{proof}

The final step in the proof of Theorem 1$^\prime$ is the next
proposition.

\begin{prop}
\label{prop:constructtransversal} Let $\cF$ be an orientable and
transversely orientable foliation of codimension 2, $B_{0} \supset
B_{1} \supset \cdots$ its Epstein hierarchy, and $N_\al$ the
union of all non-compact leaves of $\cF$ in $B_\al\setminus B_{\al+1}$. Assume
that $B_{1}$  is compact, that 
$\bigcap_\al B_\al$ is empty or consists of non-compact leaves only,  and
that 
$\bigcap_\al B_\al$ and all
$N_\al$ have dimension at most equal to dim$\cF$. Then there exists a compact
transverse 2-manifold $T$ intersecting every leaf of $B_{1}$ such that $\partial
T \cap B_{1} = \emptyset$.
\end{prop}

\begin{proof} In the absence of non-compact leaves (when $\bigcap_\al B_\al$ and
all $N_\al$ are empty) the proposition was proved in \cite{EMS} and \cite{Vo1}
by extending the key ideas of Epstein in \cite{Ep1}. Our
proof here is basically the same by noticing at each step that the non-compact
leaves cause no additional difficulties.

Assume first that $\bigcap_\al B_\al = \emptyset.$  Then by
Proposition~\ref{prop:hierarchyend} there exists an ordinal $\ga$ such that
$B_{\ga} \neq
\emptyset$ and $B_{\ga + 1} =
\emptyset$. We may assume that $\ga \geq {1}$, for otherwise there is nothing to
prove. Then $B_{\ga}$ is compact and again by \ref{prop:hierarchyend} contains
only compact leaves. Therefore $B_{\ga}$ is a Seifert fibration with an at most
1-dimensional leaf space. The techniques of \cite{EMS}  and \cite{Vo1} then
show that there exists a compact transverse manifold $T_{\ga}$
intersecting each leaf of $B_{\ga}$ such that $\partial T _{\ga} \cap
B_{\ga} = \emptyset$. For a detailed proof see \cite{Vo2}, Proposition
4.7.

If $\bigcap_\al B_\al \neq \emptyset$ there
exists an ordinal
$\gamma$ such that $B_\gamma = \bigcap_\al B_\al$. By hy\-pothe\-sis $B_\gamma$
is transversely 0--dimensional and thus we can again find a compact transverse
$T_\gamma$ with the properties above.

The idea is now to use downward induction, i.\ e., if
$\alpha > 1$, and  if $T_{\al}$ is a compact transverse 2-manifold 
which intersects every leaf of $B_{\al}$ and whose boundary
$\partial T_{\al}$ is disjoint from $B_{\al}$,  we have to
construct for some $\beta <\al$ a transverse 2-manifold
$T_{\beta}$ having the same properties with respect to $B_\beta$. If $\al$ is
a limit ordinal then for some $\beta < \al$ the 2-manifold $T_\alpha$
intersects every leaf of
$B_{\beta}$ and
$\partial T_\alpha \cap B_{\beta} = \emptyset$. This can be seen as follows.

The union $A$ of leaves of $\cF$ not intersecting int$T_\alpha$ is closed.
Therefore,  for any $\beta \geq 1$ the set $(A \cup \partial T_\alpha ) \cap
B_{\beta}$ is compact and $\bcap_{\delta < \al} ((A \cup \partial T_\alpha )
\cap B_{\delta}) = (A \cup \partial T_\alpha ) \cap B_{\al} = \emptyset$. It
follows that for some $\beta < \al$ we have $(A \cup \partial T_\alpha ) \cap
B_{\beta} = \emptyset$ which implies that $T_\alpha$ intersects every leaf of
$B_{\beta}$ and $B_{\beta} \cap \partial T_\alpha = \emptyset$.

So we may assume that $\al$ is not a limit ordinal. By
Lemma~\ref{lem:avoidboundary} we may also assume that for the union $N_{\al -
1}$ of all  non-compact leaves of $B_{\al-1}\setminus B_{\al}$ we have
$N_{\al-1}
\cap
\partial T_{\al} = \emptyset$.

Now, $B_{\al-1} \smetmi (N_{\al-1} \cup B_{\al})$ is a Seifert
fibration. Since $\partial T_{\al}$ is compact we find a closed
invariant neighborhood $K_{0}$ of $B_{\al} \cup N_{\al-1}$ in $B_{\al
-1}$ such that $T_{\al}$ intersects every leaf of $K_{0}$ and $K_{0}
\cap \partial T_{\al} = \emptyset$. Since exceptional leaves of
 foliated Seifert fibred subsets of a $C^1$-foliation are
isolated (\cite{Vo2}, Lemma 4.4) we may also assume that the set theoretic
boundary
$Fr_{B_{\al-1}} (K_{0})$ of $K_{0}$ in $B_{\al-1}$ does not contain any
exceptional leaf of the Seifert fibration $\cF \mid \big (B_{\alpha-1} \smetmi
(N_{\al-1} \cup B_{\al}) \big )$. (As a reminder: a leaf of a Seifert fibration
is called exceptional, if its holonomy group is non-trivial.)

Below we will establish the following claim.

\begin{claim}
\label{claim:exceptionalleaves} Let $E$ be the union of the exceptional leaves
of the Seifert fibration $\cF \mid (B_{\al-1} \smetmi K_{0})$. Then there
exists a compact invariant neighborhood $N$ of $E$ in
$B_{\al-1} \smetmi K_{0}$ and a compact transverse manifold $S$ with
the following properties

\bitem
\item[\rm (i)]  $S$ intersects every leaf of $K_{1} = K_{0} \cup N$;
\item[\rm (ii)] $\partial S \cap K_{1} = \emptyset$~;
\item[\rm (iii)]  there exists a $\rho > 0$ and an invariant
neighborhood $U_{1}$ in $B_{\al-1}$ of the point set theoretic boundary
$Fr_{B_{\al-1}} K_{1}$
 of $K_1$ in
$B_{\al-1}$ such that every leaf of $U_{1}$ intersects $S$ in exactly $\rho$
points. 
\eitem
\end{claim}

Assuming that \ref{claim:exceptionalleaves} is true we then proceed as in [Ep
1], [EMS], [Vo 1], [Vo 2] to extend $S$ to a transverse manifold having
properties (i) and (ii) above with $K_{1}$ replaced by $B_{\al-1}$. The idea is
to cover the locally trivial bundle $Cl (B_{\al-1} \smetmi K_{1})$ by
finitely many bundle charts $C_{2}, \ldots, C_{n}$ and then to construct
inductively transverse compact manifolds $S_1 = S, S_2, \ldots, S_n$ such that
$S_i$ has properties (i), (ii), and (iii) above with $K_1$ replaced by $K_i =
K_1 \cup C_2 \cup \cdots \cup C_i$. This is done by choosing for each $C_i$ a
compact transverse manifold $D_i$ with $\partial D_i \cap C_i = \emptyset$ and
intersecting each leaf of $C_i$ in exactly $\rho$ points.  Then we
shrink at each step $S_{i}$ and $D_{i+1}$ somewhat  and adjust $D_{i+1}$ so
that $S_{i+1} = S_{i} \cup D_{i+1}$ is a transverse 2-manifold having properties
(i), (ii), and (iii) with regard to $K_{i+1}$. For a detailed description of this
see \cite{Vo2}, proof of 4.7 (Note that in figure 2 of \cite{Vo2} each $\Gamma$
should be interpreted as the intersection symbol $\cap$). 
\end{proof}

\begin{proof}[Proof of \ref{claim:exceptionalleaves}] (An adaptation of the
proof in \cite{Ep1}, Section 10, to our situation.) Since $K_{0} \cap \partial
T_{\al} =
\emptyset$, since
$K_{0}$ is a neighborhood of $N_{\al-1} \cup B_{\al}$ in $B_{\al-1}$
and since $Fr_{B_{\al-1}} (K_{0})$ does not contain an exceptional leaf
of $B_{\al-1} \smetmi (N_{\al-1} \cup B_{\al})$, we find an invariant
compact neighborhood $V = V_{1} \cup \cdots \cup V_{k}$ of $Fr_{B_{\al-1}}
(K_{0})$ in $B_{\al-1}$ such that $V \cap E = \emptyset$, the $V_{i}$ are
disjoint compact invariant sets and $\sec_{T_{\al}}$ restricted to each
$V_{i}$ is constant with value, say $n_{i}$. We will assume that the $n_{i}$ are
pairwise distinct. Let $U = U_{1} \cup \cdots \cup U_{k}$ be another compact
invariant neighborhood of $Fr_{B_{\al-1}} (K_{0})$ such that for
all $i$ we have $U_{i} \subset \mbox{int}_{B_{\al-1}} (V_{i})$. Then
every component $C$ of $K_{0}$ which is not entirely contained in
$V_{i}$ and meets $U_{i}$ has infinitely many leaves in $V_{i}$.
Our hypotheses let us apply  Proposition~\ref{prop:Weaver}  to components
of $K_0\cap T_\al$. From this we conclude  that no component of
$K_{0}$ will intersect two of the sets $U_{i}$. It is now a routine matter
(see \cite{Ku}, \S 47 II Theorem 3) to decompose $K_{0}$ into disjoint closed
subsets $K_{0, 1} \cup \cdots \cup K_{0, k}$ such that $K_{0, i} \cap
U_{j}$ is empty for $i \neq j$. The closure of $B_{\al-1} \smetmi
K_0$ is compact, and $\cF$ restricted to this set
is a Seifert fibration. Therefore, $E$ is a finite union of leaves
$L_{1}, \ldots, L_{m}$. Let $U_{k+i}$ be a compact invariant
neighborhood of $L_{i}$ in $B_{\al-1}$. If the $U_{k+i}$ are
small enough we may assume that all leaves of $U_{k+i}$ intersect a transverse
disk $D_{i}$ in exactly $n_{k+i}$ points except $L_{i}$ which intersects
$D_{i}$ once. We may further assume that $D_{i} \cap K_{0} = \emptyset$,
$U_{k+i} \cap \partial D_{i} = \emptyset$ and that for $i \neq j ~ U_{k+i} \cap
D_{j} = \emptyset$. Now let $T_{i},\, i = 1, \ldots, k$, be a compact
2-manifold-neighborhood of $K_{0,i} \cap T_{\al}$ having the following properties

\begin{tabular}{rll}
(i) & $K_{0,i} \cap T_{j}  = \emptyset$  & ,  $i \neq j$;\\[2mm]
(ii) & $K_{0,i} \cap \partial T_{i}  =  \emptyset$ & ,  $i = 1, \ldots,
k$;\\[2mm]
(iii) & $T_{i} \cap U_{k+j}  =  \emptyset$ & ,  $i = 1, \ldots, k$; $j = 1,
\ldots, m $.  
\end{tabular}

\noindent
It is clear that we can find the $T_{i}$ with the desired properties.
Since $\cF$ is orientable and since we may assume that every
component of every $T_{i}$ is a compact 2-manifold with non-empty
boundary, tubular neighborhoods of the $T_{i}$ and $D_{i}$ are
trivial. We also may assume that all fibres of these neighborhoods are
open disks in leaves of $\cF$. We find the desired transverse
manifold $S$ by replacing each $T_{i}$ by $m_{i} = (n_{1} \cdot~
\cdots~ \cdot n_{k+m})/n_{i}$ disjoint copies of $T_{i}$ each
being a section of the tubular neighborhood of $T_{i}$ and similarly
$D_{i}$ by $(n_{1} \cdot~ \cdots~ \cdot n_{k+m})/n_{k+i}$ copies,
making sure that all these copies are disjoint. Then $S$ is the union of all
these copies and $N = U_{k+1} \cup \cdots \cup U_{k+m}$.  
\end{proof}

\noindent
{\large\bf Appendix:\qua the topological case}
\medskip

The hypothesis that $\cF$ is $C^{1}$ was used in the preceding
section in two instances. First in the proof of \ref{prop:Weaver}
and second in the statement that exceptional leaves of a Seifert fibration are
isolated. The first instance can be dealt with as in the proof of
Theorem 3.3 in \cite{Vo2}. Instead of showing that each $D_{m} \smetmi
D_{m-1}$ is closed one shows that any component $R$ of $D_{m}$ which
meets $D_{m-1}$ is entirely contained in $D_{m-1}$. This is Lemma 3.5
in \cite{Vo2} and its proof can be used in our situation since
non-compact leaves do not figure in $D_{m}$.

Exceptional leaves need not be isolated in topological Seifert
fibrations. See Remark 4.5 in  \cite{Vo2}. There are two ways to get
around this problem. One is to show that nevertheless we can argue as
before in the proof of Claim~\ref{claim:exceptionalleaves} by showing the
existence of compact invariant neighborhoods $U$ of these leaves with the
following property: there is an invariant neighborhood $F$ of
$Fr_{B_{\al-1}} U$ in $B_{\al -1}$  such that all leaves in $F$ will intersect a
transverse 2-manifold
$D$ in the same number of points. Here $D$ is supposed
to meet every leaf of $U$ and $U \cap \partial D = \emptyset$. The proof of the
existence of  such a $U$ follows from the arguments at the beginning of the
proof of
\ref{claim:exceptionalleaves} where we decomposed
$K_{0}$  into $K_{0,1} \cup
\cdots \cup K_{0, k}$. Decompose $U$ by the same process into $U_1\cup\cdots
U_r$ and then replace $U$ by the $U_i$  containing the exceptional leaf. This
will have the required properties.

Another way to proceed is the use of the so called fine Epstein
hierarchy instead of our version. Here $B_{\al + 1}$ is defined to be
the union of leaves $L$ in $B_{\al}$ such that for any open
transverse manifold $T$ intersecting $L$ there are leaves of $B_{\al}$
intersecting $T$ in more than one point. Then $\cF$ restricted to
$B_{\al} \smetmi (N_{\al} \cup B_{\al + 1})$ will be a locally
trivial bundle and the problem of exceptional leaves disappears
altogether.

\section{The foliation cycle and the
proof of Theorem $\mathbf{2'}$}
\label{sec:codimp}

We begin with a proof of the analogue of what is called the ``Moving
Leaf Proposition'' in \cite{EMS}.

\begin{prop}
\label{prop:movingleaf} Let $\cF$ be a foliation of codimension $k$
on a manifold $M$, let $B_{1}$ be the first bad set of $\cF$ and let
$N_{0}$ be the union of all non--compact leaves of $\cF$ in the complement
of $B_{1}$. Assume that $N_{0} \cup B_{1}$ is not empty, that the
closure of every leaf of $\cF$ is compact and that one of
the following two conditions holds

\begin{itemize}
\item[\rm(i)] $N_{0}$ is not open, or
\item[\rm(ii)] $N_{0} = \emptyset$, $\mathrm{int}B_1 = \emptyset$, and
$B_1\setminus B_2
\neq \emptyset$. 
\end{itemize}

Then for any transverse $k$--manifold $T$ whose interior intersects every leaf
of $B_{1}$ there exists a component $U$ of $M \smetmi
(N_{0} \cup B_{1})$ such that $\sec_{T}$ is unbounded on $U$ (for
notations see \ref{not:transverse} and \ref{not:sec}).
\end{prop}

\begin{proof} Assume that (i) holds. Since $B_{1}$ is closed we find
$x_{0} \in N_{0} \cap Fr N_{0}$, where $Fr N_{0}$ is the set
theoretic boundary of $N_{0}$, and a neighborhood $V$ of $x_{0}$ with
$V \cap B_{1} = \emptyset$. Let $y_{0}$ be a point of $V \smetmi N_{0}$ and
$U$ be the component of $M \smetmi (N_{0} \cup B_{1})$ containing
$y_{0}$. Then $\sec_{T}$ will be unbounded on U for any transverse
$k$--manifold $T$ such that $\mathrm{int}T$ intersects every leaf of $B_{1}$.
To see this let $x_{1}$ be a point in $(Fr U) \cap V$. Since $N_{0}
\cup B_{1}$ is closed $x_{1}$ lies in $N_{0}$. Since the closure of
the leaf $L_{x_{1}}$ through $x_{1}$ is compact (by hypothesis) there
exists a leaf $L \subset B_{1}$ in the limit set of $L_{x_{1}}$.
Therefore $L_{x_{1}}$ will intersect any transverse manifold $T$ in
infinitely many points if Int $T \cap L \neq \emptyset$. Since $L_{x_{1}}
\subset Fr U$ the function $\sec_{T}$ will be unbounded on $U$ by
Property~\ref{proper:sec}(a) of $\sec_{T}$.

If (ii) holds we distinguish two cases.

\noindent
{\bf Case 1}\qua All leaves of $B_{1} \smetmi B_{2}$ are non--compact.
Since int$B_{1} = \emptyset$, the union $N_{1} = B_{1} \smetmi B_{2}$ of all 
non--compact leaves in the complement of $B_2$ is not open unless it is empty.
By hypothesis, $B_{1} = B_{1} \cup N_{0}$  and $N_{1} = B_{1} \smetmi B_{2}$ are
not empty. Now, we can argue as before,  replacing $N_{0}$ by $N_{1}$, and
$B_{1}$ by $B_{2}$. In this way we find a component $U$ of $M \smetmi (N_{1}
\cup B_{2}) = M
\smetmi (N_{0} \cup B_{1})$ which has the following property:  $\sec_{T}$
is unbounded on $U$ for any transverse $k$--manifold $T$ whose interior
intersects every leaf of
$B_{2}$.

\noindent
{\bf Case 2}\qua $B_{1} \smetmi B_{2}$ contains compact leaves. Since
the union $N_{1}$ of  all non--compact leaves of $B_{1} \smetmi B_{2}$ is a
closed subset of $M \smetmi B_{2}$ the space $M \smetmi (N_{1} \cup
B_{2})$ is a manifold and the restriction $\cF_{1}$ of $\cF$ to $M
\smetmi (N_{1} \cup B_{2})$ is a foliation with all leaves compact.
Furthermore, the first bad set of $\cF_{1}$ is $B_{1} \smetmi (B_{2}
\cup N_{1})$ and therefore not empty.

The Moving Leaf Proposition in \cite{EMS} requires the bad set to be
compact and $B_{1} \smetmi (B_{2} \cup N_{1})$ need not be compact.
Now there are two parts in the proof of the Moving Leaf Proposition in
\cite{EMS}. The first (and most difficult) part states that there is a component
of the complement of the first bad set on which the volume of leaf
function is not bounded. The proof does not make any use of the
compactness of $B_{1} (\cF_{1})$. It is purely local. In fact what is
proved in \cite{EMS} in the two paragraphs starting with the last
paragraph on page 23 can be stated as follows: Let $\cG$ be a foliation of
codimension $k$  with all leaves compact and $L$ any leaf in $B_{1} (\cG)$ such
that
$L$ has trivial holonomy in the foliated set $B_{1} (\cG)$. Let $D$
be any transverse $k$--disk intersecting $L$ in its interior
$\overset\circ D$. Then there exists a 
component $V$ of $\overset\circ D \smetmi B_{1} 
(\cG)$ such that $\sec_{D}$ is unbounded on $V$.

Since the union of leaves of $B_{1} (\cG)$ with trivial holonomy in
$B_{1} (\cG)$ is open and dense, our claim is an immediate consequence
of the above statement when applied to $\cG = \cF_{1}$. 
\end{proof}

Next we will construct a particular foliation cycle. This is the
point where compactness of $N_{0} \cup B_{1}$ is essential.
Compactness of $N_{0} \cup B_{1}$ guarantees the existence of
arbitrarily small saturated compact neighborhoods $X$ of $N_{0} \cup
B_{1}$. This is due to the fact that the frontier $Fr (W)$ of a
relatively compact neighborhood $W$ of $N_{0} \cup B_{1}$ is a
compact subset of $M \smetmi (N_{0} \cup B_{1})$, and on $M \smetmi
(N_{0} \cup B_{1})$ the foliation $\cF$ is a Seifert fibration.
Therefore, the saturation $S$ of $Fr (W)$ is also compact and thus
closed. Then $Y = W \smetmi S$ is a saturated neighborhood of $N_{0}
\cup B_{1}$ with $Cl (Y) \subset Cl (W)$.

From now on we will assume that $N_{0} \cup B_{1}$ is compact and non-empty, that
$N_{0}$ is either empty or not open und that $\bigcap_\al B_\al = \emptyset$.
Furthermore, we assume that $\cF$ is $C^{1}$ and oriented. The last condition
allows us to consider the compact leaves of $\cF$ as $(n-k)$--dimensional cycles.

Let $X$ be a compact saturated neighborhood of $N_{0} \cup B_{1}$.
Then we can find finitely many foliation charts $W_{i} = E_{i} \times
T_{i},\; i = 1, \ldots, s$, whose interiors cover $X$. Here we assume
that each $T_{i}$ is a compact transverse $k$--manifold and each
$E_{i} \times \{t\}$ is an open relatively compact subset of a leaf.
As usual, we assume that each $ E_{i} \times T_{i}$ is part of a
larger foliation chart $\tilde{E}_i \times \tilde{T}_{i}$ with
$T_{i} \subset$ int $\tilde{T}_{i}$ and $Cl (E_{i}) \subset$ int
$\tilde{E}_{i}$. We may and will assume that the $T_{i}$ are
disjoint. Then $T = \bigcup T_{i}$ is a compact transverse $k$--manifold whose
interior int $T = \bigcup$ int $T_{i}$ intersects every leaf of $B_{1}$.

Let $U$ be a component of $X \smetmi (B_{1} \cup N_{0})$ such that
$\sec_{T}$ is unbounded on $U$. By Propositions~\ref{prop:movingleaf} and
\ref{prop:openbadset} such a component exists. Let $L_{1}, L_{2}, \ldots$ be a
sequence of leaves in
$U$ such that $\sec_{T} (L_{i})$ is a strictly increasing  unbounded sequence. 
Since $\sec_{T}$ is bounded on any compact subset of $U$ the sequence of
leaves $L_{1}, L_{2}, \ldots$ converges to $B_{1} \cup N_{0}$. Since
the union of all leaves of $U$ with trivial holonomy is open and dense we
may and will assume that the leaves $L_{i}$ have trivial holonomy.
Then all leaves $L_{i}$ are homologous in $U$.

In \S\S \,2 and 3 of \cite{EMS} is explained how this set--up leads after
passing to a suitable subsequence and the appropriate choice of integers
$n_{i}$ to a limiting foliation cycle $\displaystyle \lim_{i}
\frac{1}{n_{i}} L_{i}.$ This foliation cycle will be essential in the
proof of Theorem $2'$. We repeat its construction. For each $i$ let $n_{i} = \max
\{\sec_{T_{j}} (L_{i}) : j = 1, \ldots, s\}$. By passing to a subsequence of the
$L_{i}$ and reordering the $T_{j}$ we may assume that $n_{i} = \sec_{T_{1}}
(L_{i})$. We define a non--negative measure $\mu_{j, i}$ on the Borel sets of
$T_{j}$ by assigning each point of $T_{j} \cap L_{i}$ the mass
$\frac{1}{n_{i}}$. Then $\mu_{j, i} (T_{j}) \leq 1$ for all $i, j$ and $\mu_{1,
i} (T_{1}) = 1$ for all $i$. Consequently, after passing to a further subsequence
of the $L_{i}$, we may assume that for all $j$ the measures $\mu_{j, i}$ converge
to a non--negative measure $\mu_{j}$ on $T_{j}$ with $\mu_{j} (T_{j}) \leq 1$
and $\mu_{1} (T_{1}) = 1$.

By Lemma A of \S \,3 of \cite{EMS} the measures $\{\mu_{j}\}$ are holonomy
invariant and therefore define a geometric current $C 
\{\mu_{j}\}$. The associated closed de~Rham current is equal to
$\displaystyle \lim_{i}\frac{1}{n_{i}} L_{i}$, i.\ e.\ 
for any
$(n-k)$--form $\om$ defined in a neighborhood of $X$ we have $\lgl
C \{ \mu_{j} \}, \om \rgl = \displaystyle \lim_{i} \frac{1}{n_{i}}
\int_{L_{i}} \om$. This is Lemma B of \cite{EMS}.  From this we obtain the first
important property of our foliation  cycle $C \{\mu_{j}\}$.
\vspace{2mm}

\begin{noname}[Property 1 of the foliation cycle $C \{\mu_{j} \}$]
 Let $\om$ be any closed $(n-k)$--form defined in $\mathrm{int}X$ then
$\lgl C \{\mu_{j} \}, \om \rgl = 0$.
\end{noname}

\begin{proof} This is due to the simple fact that the leaves $L_{1},
L_{2}, \ldots$ are all homologuous in $\mathrm{int}U \subset\mathrm{int}X$ so
that the sequence $\int_{L_{i}} \om$ is constant if $\om$ is closed. Since
$1/{n_{i}}$ converges to $0$ we are done.  
\end{proof}

The proof of Theorem $2'$ is now an immediate consequence of the second
property of $C \{\mu_{j}\}$.

\begin{noname}[Property 2 of the foliation cycle $C \{\mu_{j} \}$]
\label{cycleprop2}
 Let $\om$ be any closed $(n-k)$--form defined in a neighborhood of
$X$ such that for any compact leaf $L$ of $B_{1}$ the inequality
$\int_{L} \om > 0$ holds. Then
$$\lgl C \{\mu_{j} \}, \om \rgl > 0~.$$
\end{noname}

\begin{proof} Recall the definition of the de Rham current 
$$
\lgl C\{\mu_{j} \},-\rgl:\linebreak \Om^{n-k}(M_{0})~\lra~\R,
$$ 
where $M_{0}$ is any neighborhood of $X$. One chooses a partition of unity
$p_{1}, \ldots, p_{s}$ subordinate to the covering of $X$ by the interiors of
$W_{j} = E_{j}
\times T_{j}$ and defines for any $\eta \in \Om^{n-k} (M_{0})$
$$\lgl C \{\mu_{j}\}, \eta \rgl = \Sum_{j} \Int_{T_{j}} \left (
\Int_{E_{j} \times \{t\}} (p_{j}\cdot\eta) \right ) d{\mu_{j}}
(t)~.$$
The definition is easily seen to be independent of the choice of 
partition of unity. \cite{EMS}, \S\,2.

We need to change the ``local'' recipe for calculating $\lgl C
\{\mu_{j}\}, \eta \rgl$ to a more global one where we integrate
$(n-k)$--forms over total leaves instead of plaques $E_{j} \times
\{t\}$.

First we notice that for every $j$ the measure $\mu_{j}$ is supported
on $T_{j} \cap B_{1}$. Clearly, $\mu_{j}$ is supported on $(N_{0} \cup
B_{1}) \cap T_{j}$ since the sequence $L_{i}$ converges to the closed
set $N_{0} \cup B_{1}$. Let $x$ be a point of $N_{0} \cap T_{j}$ and
$\tilde{T}_{j}$ a transverse $k$--manifold such that $T_{j} \subset$
Int $\tilde{T}_{j}$. Then we find a transverse $k$--manifold $D
\subset\mathrm{int}\tilde{T}_{j}$ such that $x \in\mathrm{int}D$ and $\sec_{D}$
is bounded. In particular, the number of intersection points of $L_{i}$
with $D$ is bounded, and this implies that $\mu_{j} (D) = 0$.

By hypothesis $\bigcap_\al B_\al = \emptyset$. Then
Proposition~\ref{prop:hierarchyend} tells us that
$B_{1}$ is a {\it countable} disjoint union of Borel sets:
\begin{center}
$B_{1} = \bcup_{\al \geq 1} (B_{\al} \smetmi B_{\al + 1})~.$
\end{center}
As before, denote the union of all  non--compact leaves of $B_{\al} \smetmi
B_{\al + 1}$ by $N_{\al}$. Then we make the following claim.

\begin{claim} For all $\al$ and $j$ the equation $\mu_{j} (N_{\al}
\cap T_{j}) = 0$ holds.
\end{claim}

\begin{proof}[Proof of Claim] We have already proved this statement for $\al =
0$ using an easy argument. In outline, the statement is true in general
because leaves of $N_{\al}$ have their limit points in $B_{\al+1}$,
and, if $B_{\al+1} \cup N_{\al}$ is
compact, their limit sets are non--empty. If $\mu_{j} (N_{\al} \cap T_{j}) > 0$
for some $j$, we find a compact  set $E\subset N_\al \cap T_j$ in the complement
of
$B_{\al+1}$ with  $\mu_{j} (E) > 0$.  Using holonomy translations repeatedly we
can push
$E$ into a countable disjoint family of subsets in
$\bigcup_i T_i$. The holonomy invariance of the measures then implies that each
of these sets has measure not less than $\mu_{j} (E)$.
 This will
contradict the fact that for all $i$ we have $\mu_{i} (T_i) \leq 1$.

In more detail, assume that $\mu_{j} (N_{\al} \cap T_{j}) > 0$. Then by
passing, if necessary, to a different $T_j$ we may also assume that $\mu_{j}
(N_{\al} \cap\mathrm{int}T_{j}) > 0$. The set
$N_{\al} \cap\mathrm{int}T_{j}$ is covered by (countably many) sets of the form
$N_{\al} \cap S$ such that $S$ is open in $\mathrm{int}T_j$ and
$\sec_{S}$ is bounded on
$N_{\al}$. So
we may assume that $\mu_{j} (S \cap K_0) > 0$ for some such $S$ and
some compact subset $K_0$ of the closed subset $S \cap N_{\al}$ of $S$.
By (\ref{proper:sec}) $\sec_{S}$ is lower semicontinuous and thus by the Baire
category theorem there exists an open dense subset of $K_0$ where
$\sec_{S}$ is continuous. Let $K_{1}$ be its complement in $K_0$. Then
clearly $\max \{\sec_{S} (x) \mid x \in K_{1}\} < \max \{\sec_{S} (x)
\mid x \in K_0\}$. Continuing inductively we find a finite sequence $K_0
\supset K_{1} \supset K_{2} \supset \cdots \supset K_{r} =
\varnothing$ such that $\sec_{S}$ is locally constant on $K_{i}
\smallsetminus K_{i+1}$ for all $i \geq 0$. Using
this, we find an open subset $U$ of $T_{j}$ and a compact subset $N$
of $U \cap N_{\al}$ with $\mu_{j} (N) > 0$ such that $\sec_{U} (x) = 1$ for
all $x \in N$.

Now, each leaf through a point of $N$ accumulates against
$B_{\al+1}$. Using holonomy translations along the leaf through $x
\in N$ we find a curve $\omega$ through this leaf and a neighborhood $V
(x)$ of $x$ in $N$ such that $\omega$ holonomy--translates $V (x)$
into a set, $V' (x)$, contained in some $T_{i}$ such that $N \cap V' (x)
= \varnothing$. By compactness of $N$ finitely many of such $V (x)$
suffice to cover $N$. Assuming that each $V (x)$ is compact the union of these
$V' (x)$ constitute a compact set $N'$ disjoint from $N$ such that 
$\sum\limits_{i} \mu_{i} (N' \cap T_{i})
\geq \mu_{j} (N \cap T_{j})$. The inequality is due to the fact
that $\sum \mu_{i}$ is holonomy invariant and that each leaf of
$N_{\al}$ intersects $N$ in at most one point. Therefore, if $V'
(x_{1})$ and $V'(x_{2})$ intersect, then $V(x_{1})$ and $V(x_{2})$
intersect in a set of at least the same measure. The inequality then follows by
induction on the number of $V(x_i)$ used to cover $N$.  Since
$N'$ is compact, we can do the process over again, moving open sets $W (x)$ of
$N$  to sets $V'' (x)$ in
some
$T_{i}$ such that $V'' (x) \cap (N \cup N') = \varnothing$, obtaining 
a set $N''$ such that $\sum\limits_{i} \mu_{i} (N'' \cap T_{i}) \geq
\mu_{j} (N)$. Continuing, we find compact subsets of
$\bigcup\limits_{i} (N_{\al} \cap T_{i})$ of arbitrary large measures
contradicting $\mu_{i} (T_{i}) \leq 1$ for all $i$. 
\end{proof}

Therefore, for any $j$
\begin{center}
$\mu_{j} (T_{j}) = \Sum_{\al \geq 1} \mu_{j} [\big ( B_{\al} \smetmi
(B_{\al + 1} \cup N_{\al}) \big ) \cap T_{j}]~.$
\end{center}
\noindent
If $\mu_{\al, j}$ denotes the restriction of $\mu_{j}$ to $\big
(B_{\al} \smetmi (B_{\al+1} \cup N_{\al}) \big ) \cap T_{j}$, then
$\{\mu_{\al, j} \}, j = 1, \ldots, s$,
defines for any $\al \geq 1$ a holonomy invariant transverse measure
on $T = \bcup_{j} T_{j}$. We denote by $C \{\mu_{\al, j}\}$ the
associated foliation cycle. Then
\begin{center}
$C \{\mu_{j}\} = \Sum_{\al \geq 1} C \{\mu_{\al, j}\}~.$
\end{center}
The proof of (\ref{cycleprop2}) is thus a consequence of the next lemma.
\end{proof}

\begin{lem}
\label{lem:positivemeasure} Let $\omega$ be any $(n-k$)--form defined in a
neighborhood of $B_{1}$ such that $\int_{L} \omega > 0$ for any
compact leaf of $B_{1}$. Then for any $\al \geq 1$ we have $\lgl C
\{\mu_{\al,j}\}, \omega \rgl \geq 0$ and there exists at least one
$\al \geq 1$ such that $\lgl C \{\mu_{\al, j}\}, \omega \rgl > 0$.
\end{lem}

\begin{proof} A proof of this lemma can easily be extracted from
\S\S\ 6 and 7 in \cite{EMS}. Our set--up is slightly different. In
particular, we use what in \cite{EMS} is called the coarse Epstein
filtration. For the convenience of the reader we give a direct proof
of \ref{lem:positivemeasure} adjusted to our situation. Fix $\al \geq 1$. The
foliation
$\cF$ restricted to $S_{\al} = B_{\al} \smetmi (B_{\al+1} \cup N_{\al})$ is a
foliation with all leaves compact. By definition  of the Epstein hierarchy (see
\ref{defn:hierarchy}) there exists for any leaf
$L$ of $S_{\al}$  a transverse disk $D$ such that $L$
intersects \ int$D$ and $\sec_{D}$ is bounded on $S_{\al}$. This implies
that $\cF$ restricted to $S_{\al}$ is a Seifert fibration which
translates into a very explicit description of a foliated
neighborhood $U_{L}$ of $L$ in $S_{\al}$ as follows (for details see
\cite{Ep2}). Let
$D$ be a transverse manifold whose interior intersects $L$. Let  $x$ be a
point  in 
$\mathrm{int} D \cap L$. Then we find a neighborhood $U_{x}$ of $x$ in $S_{\al} \cap
D$, a finite group $H$ of homeomorphisms of $U_{x}$ fixing $x$, and a finite regular
covering
$\tilde{L} \lra L$ with deck transformation group isomorphic to $H$
such that as a foliated set $U_{L}$ is isomorphic to $(\tilde{L}
\times U_{x})/H$. Here $H$ operates diagonally on $\tilde{L} \times
U_{x}$ and $(\tilde{L} \times U_{x})/H$ is foliated by the images of
$\tilde{L} \times \{t\}, t \in U_{x}$. By choosing $U_{x}$
sufficiently small we may assume that for any $h \in H$ the germ of
$h$ at $x$ is not trivial, so that $H$ realizes the holonomy group of
$L$ in the foliated set $S_{\al}$.

Fixing again a leaf $L$ of $S_{\al}$ we find an index $j_{0}$ such
that $L \cap$ int$T_{j_{0}} \neq \emptyset$. Then we may choose $x \in L
\cap$ int$T_{j_{0}}$ and $U_{x} \subset$ int$T_{j_{0}}$ in the
above discussion to have the following additional properties:

\begin{proper}
\label{proper:holtranslate} For any $j \in \{1, \ldots, 
s\}$ and $y \in L \cap\mathrm{int}T_j$ there exists a holonomy translation
$$h_{y x} : U_{x} \lra \mathrm{int}T_{j}$$
along a path in $L$ from $x$ to $y$ such that
\bitem
\item[\rm (i)]  $h_{yx} (U_{x}) \cap h_{y' x} (U_{x}) = \emptyset$,\; if $y \neq
y'$; 
\item[\rm (ii)] Let $\{p_{j}\}$
be the partition of unity subordinate to $\{\mathrm{int}W_{j} =
E_{j}\times\mathrm{int}T_{j}\}$ used in our formula for evaluating our
foliation cycle on forms. Then for every $j$ the projection to
$\mathrm{int}T_{j}$ of the intersection of the neighborhood $U_{L} =
(\tilde{L}\times~U_{x})/H$ of
$L$ with 
$\mathrm{supp}(p_{j}) \subset E_{j} \times \mathrm{int}T_{j}$ is
contained in

$\bcup_{y \in \mathrm{int} T_{j}} h_{yx} (U_{x})~.$
\eitem
\end{proper}

Using these properties we can rewrite the contribution of $U_{L}$ to
$C \{\mu_{\al, j}\}$:

\hsp{2.4cm} 
$\Sum_{j}~~ \Int_{T_{j} \cap U_{L}} \left ( \Int_{E_{j}
\times \{t\}} p_{j} \cdot \omega \right ) d\mu_{\al, j}(t) =$

\hsp{2cm} 
$= \Sum_{j}~~ \Sum_{y \in L \cap \mbox{\scriptsize{int}}
T_{j}}~~ \Int_{U_{x}} \left ( \Int_{E_{j} \times \{h_{yx} (t) \}} p_{j}
\cdot \omega \right ) d \mu_{\al, j_{0}}(t) =$

\hsp{2cm} $= \Int_{U_{x}} \left ( \Sum_{j}~~ \Sum_{y \in L 
\cap \mbox{\scriptsize{int}} T_{j}}~~ \Int_{E_{j} \times \{h_{yx} (t)
\}} p_{j} \cdot \omega \right ) d \mu_{\al, j_{0}}(t)$.

The first equation is a consequence of the holonomy invariance of
$\{\mu_{\al, j} \}$ and property~\ref{proper:holtranslate}. The second
equation is obvious.

Also by holonomy invariance the last expression does not change if $t
\in U_{x}$ is replaced by $h (t)$ with $h \in H$. Therefore, it is
equal to
\begin{center}
$\frac{1}{\mid H \mid} ~\Int_{U_{x}} \left ( \Sum_{j}~~ \Sum_{y \in L 
\cap \mbox{\scriptsize{int}} T_{j}}~~ \Sum_{h \in H}~~ \Int_{E_{j} 
\times \{h_{yx} \circ h (t) \}} p_{j} \cdot \omega \right ) 
d \mu_{\al, j_{0}} (t)~.$
\end{center}
\noindent
Using once more (\ref{proper:holtranslate})  we see that for any $t \in U_{x}$
\begin{center}
$\Sum_{y \in L 
\cap \mbox{\scriptsize{int}} T_{j}}~~ \Sum_{h \in H}~~ \Int_{E_{j} 
\times \{h_{yx} \circ h (t) \}} p_{j}\cdot \omega = \mid H_{t} \mid \cdot~
\Int_{L_{t}} p_{j}\cdot \omega$
\end{center}
\noindent
where $L_{t}$ denotes the leaf through $t$, and $H_{t}$ is the
stabilizer of $t$ in $H$. Altogether we see that the contribution of
$U_{L} \subset S_{\al}$ to $C \{\mu_{\al, j}\}$ when evaluated on the
form $\omega$ equals

\begin{noname}
\label{formula}
\hsp{3.5cm}
$\displaystyle \frac{1}{\mid H \mid} \Int_{U_{x}} \left (
\mid H_{t}  \mid \Int_{L_{t}} \omega \right ) d \mu_{\al, j_{0}(t)}.$
\end{noname}

\noindent
Exactly the same formula holds when $U_{x}$ is replaced by an
$H$--invariant measurable subset $\tilde{U}$ of $U_{x}$ and $U_{L}$ by the union
$\tilde{U}_{L}$ of leaves through $\tilde{U}$. Now, $S_{\al}$ is a
countable disjoint union of sets of the form $\tilde{U}_{L}$. To see
this cover $S_{\al}$ by a locally finite (and therefore countable)
family of open Seifert fibred neighborhoods $U_{L_{k}}$ of leaves
$L_{k}$ having all the properties needed for the discussion above and
having compact closure in $S_{\al}$. Then $S_{\al}$ is the disjoint
union of $U_{L_{k}} \smetmi \bcup_{i < k} U_{L_{i}}, k = 1, 2,
\ldots$~. 

The proof of Lemma~\ref{lem:positivemeasure} is now immediate from
(\ref{formula}). The hypothesis that for any compact leaf $L \subset B_{1}$ the
integral
$\Int_{L} \omega$ is positive guarantees that the expression in (\ref{formula})
is never negative. On the other hand, for some\, $\al\geq 1$\, and some\,
$j_{0}$\, the measure\, $\mu_{\al, j_{0}}\, ({\rm int\,} T_{j_{0}})\, =\,
\mu_{j_{0}}\,
({\rm int\,} T_{j_0} \cap S_{\al})$ is positive since $\{\mu_{j}\}$ is
holonomy invariant, $\mu_{1}\, (T_{1}) = 1$ and $\{\mathrm{int}T_j\times E_j\}$
covers all of $B_{1} \cup N_{0}$. Finally $t \longmapsto \mid H_{t}
\mid \Int_{L_{t}} \omega$ is continuous on $U_{x}$ and everywhere
positive. 
\end{proof}

{\makeatletter

\@thebibliography@{EMS}\small\parskip0pt % 
plus2pt\relax

\makeatother

\bibitem[EMS]{EMS}{\bf R Edwards}, {\bf K\,C Millett}, {\bf D Sullivan},
{\it Foliations with all Leaves Compact}, Topology 16 (1977) 13--32

\bibitem[Ep 1]{Ep1} {\bf D\,B\,A Epstein}, {\it Periodic flows on
3-manifolds}, Annals of Math. 95 (1972) 58--82

\bibitem[Ep 2]{Ep2} {\bf D\,B\,A Epstein}, {\it Foliations with all
leaves com\-pact}, An\-nales de l'Insti\-tut Fourier Grenoble 26 (1976) 
265--282

\bibitem[Ep 3]{Ep3} {\bf D\,B\,A Epstein},  {\it Pointwise Periodic
Homeomorphisms}, Proc. London Math. Soc., 3. Ser.42 (1981) 415--460
  
\bibitem[EMT]{EMT} {\bf D\,B\,A Epstein}, {\bf K\,C Millett}, {\bf D Tischler},
{\it Leaves without holonomy}, J.~London Math. Soc.~(2) 16 (1977)
548--552

\bibitem[Hae]{Hae} {\bf A Haefliger}, {\it
Vari\'{e}t\'{e}s feuillet\'{e}es}. Ann. Scuola Norm. Sup. Pisa 16 (1962)
367--397 

\bibitem[Ku]{Ku} {\bf K Kuratowski}, {\it Topology Vol 2}  Academic
Press 1968

\bibitem[L]{L} {\bf R Langevin}, {\it A list of questions
about  foliations}, in: {\it Differen\-tial Topology, Foliations, and Group
Actions},   P Schweitzer et al (editors), Contemporary Mathematics 161
Amer.  Math. Soc.  (1994) 59--80

\bibitem[M]{M} {\bf D Montgomery}, {\it
Pointwise periodic homeomorphisms}, Amer. J. Math. 59 (1937)
118--120

\bibitem[R]{R} {\bf G Reeb},  {\it Sur certaines pro\-pri\'{e}t\'{e}s  
topologiques des vari\'{e}t\'{e}s feuille\-t\'{e}es}, Act. Sci. Indust.,  
Hermann  1952

\bibitem[Vo 1]{Vo1} {\bf E Vogt}, {\it Foliations of Codimension 2 with
all Leaves Compact}, manu\-scrip\-ta math. 18 (1976) 187--212

\bibitem[Vo 2]{Vo2} {\bf E Vogt}, {\it Bad sets of compact
foliations of codimension~2}, Proceedings of Low Dimensional Topology
(K.~Johannson ed.) Knoxville 1992, International Press Company (1994)
187--216

\bibitem[Vo 3]{Vo3} {\bf E Vogt}, {\it Negative Euler
characteristic as an obstruction to the existence of periodic flows on open
3--manifolds}, in: T Mizutani et al (eds.), Proceedings of the international
symposium/workshop on {\it ``Geometric Study of Foliations''}, 
 World Scientific (1994) 457--474

\bibitem[Wea]{Wea} {\bf N Weaver}, {\it  Pointwise periodic homeomorphisms},
Annals of Math. 95 (1972) 83--85

\endthebibliography}

\Addresses\recd
\end{document}

%% file: agtout.tex
\def\ifplaintex{\expandafter\ifx\csname documentclass\endcsname\relax}

\def\gtp{{\mathsurround=0pt\it $\cal G\mskip-2mu$eometry \&\ 
$\cal T\!\!$opology $\cal P\!$ublications}}  % GT publications

\def\recd{{\small Received:\qua\receiveddate\ifx\reviseddate\relax
\else\qquad Revised:\qua\reviseddate\fi\par}} 

\def\lognumber#1{\def\thelognumber{#1}}
\def\volumenumber#1{\def\thevolumenumber{#1}}
\def\volumeyear#1{\def\thevolumeyear{#1}}
\def\papernumber#1{\def\thepapernumber{#1}}
\def\pagenumbers#1#2{\def\startpage{#1}\def\finishpage{#2}}
\def\published#1{\def\publishdate{#1}}

\def\received#1{\def\receiveddate{#1}}
\def\revised#1{\def\reviseddate{#1}}
\def\accepted#1{\def\accepteddate{#1}}

\def\asciiaddress#1{\def\theasciiaddress{#1}}

\long\def\asciiabstract#1{\long\def\theasciiabstract{#1}}

\let\\\par\let\thelognumber\relax\let\thevolumenumber\relax
\let\thepapernumber\relax\let\thevolumeyear\relax\let\startpage\relax
\let\finishpage\relax\let\publishdate\relax\let\receiveddate\relax
\let\reviseddate\relax\let\accepteddate\relax\let\theasciititle\relax
\let\theasciiauthors\relax\let\theasciiaddress\relax
\let\theasciiabstract\relax

\let\theasciiemail\relax

\ifplaintex
\font\logobig=cmssbx10 scaled 3836
\font\logomed=cmssbx10 scaled 2557
\else
\font\logobig=cmssbx10 scaled 4200
\font\logomed=cmssbx10 scaled 2800
\fi

\long\def\makeagttitle{   %%% start of definition of \makeagttitle
\count0=\startpage
\agt\hfill      %   Journal title (top left) 
\hbox to 45truept{\vbox to 0pt{\vglue -13truept{\logomed A\kern -.37em{\logobig 
T}\kern -.38em G}\vss}\hss}
\break
{\small Volume \thevolumenumber\ (\thevolumeyear)
\startpage--\finishpage\nl
Published: \publishdate}

\vglue .25truein

{\parskip=0pt\leftskip 0pt plus
1fil\def\\{\par\smallskip}{\Large\bf\thetitle}\par\medskip} \vglue
0.05truein

{\parskip=0pt\leftskip 0pt plus 1fil\def\\{\par}{\sc\theauthors}
\par\medskip}%
 
\vglue 0.03truein 

{\small\leftskip 25truept\rightskip 25truept{\bf Abstract}\stdspace\theabstract

{\bf AMS Classification}\stdspace\theprimaryclass
\ifx\thesecondaryclass\relax\else; \thesecondaryclass\fi\par
{\bf Keywords}\stdspace \thekeywords\par}\vglue 7truept

}   %%%% end of definition of \makeagttitle

\ifplaintex
\font\phead=cmsl9 scaled 950
\font\pnum=cmbx10 scaled 913
\font\pfoot=cmsl9 scaled 950
\headline{\vbox to 0pt{\vskip -4.5mm\line{\small\phead\ifnum
\count0=\startpage ISSN 1472-2739 (on-line) 1472-2747 (printed)
\hfill {\pnum\folio}\else\ifodd\count0\def\\{ }% 
\ifx\theshorttitle\relax\thetitle\else\theshorttitle\fi\hfill{\pnum\folio}
\else\def\\{ and }{\pnum\folio}\hfill\ifx\theshortauthors\relax\theauthors
\else\theshortauthors\fi\fi\fi}\vss}}
\footline{\vbox to 0pt{\vglue 0mm\line{\small\pfoot\ifnum\count0=\startpage
\copyright\ \gtp\hfill\else
\agt, Volume \thevolumenumber\ (\thevolumeyear)\hfill\fi}\vss}}
\else
\font=cmsl9 scaled 1050
\font\lnum=cmbx10 
\font=cmsl9 scaled 1050
\makeatletter
\def\@oddhead{{\small\lhead\ifnum\count0=\startpage ISSN 1472-2739 
(on-line) 1472-2747 (printed)\hfill {\lnum\number\count0}\else\ifodd\count0
\def\\{ }\ifx\theshorttitle\relax \thetitle \else\theshorttitle\fi\hfill
{\lnum\number\count0}\else\def\\{ and }{\lnum\number\count0}
\hfill\ifx\theshortauthors\relax 
\theauthors\else\theshortauthors\fi\fi\fi}}\def\@evenhead{\@oddhead}
\def\@oddfoot{\small\lfoot\ifnum\count0=\startpage\copyright\ \gtp\hfill\else
\agt, Volume \thevolumenumber\ (\thevolumeyear)\hfill\fi}
\def\@evenfoot{\@oddfoot}
\makeatother
\fi
\let\maketitlepage\makeagttitle
\let\makeshorttitle\maketitlepage
\let\maketitle\maketitlepage

\newwrite\gtoutfile
\long\gdef\makeheadfile{  %%% start of definition of \makeheadfile
{\def\\{, }\def\s{ }
\immediate\openout\gtoutfile head.xxx
\immediate\write\gtoutfile{To: math@arxiv.org}
\immediate\write\gtoutfile{Subject: put OR rep NNNNN:ppppp}
\immediate\write\gtoutfile{--text follows this line--}
\immediate\write\gtoutfile{Proxy-for: \ifx\theasciiauthors\relax
\theauthors\else\theasciiauthors\fi\s<\ifx\theasciiemail\relax\theemail\else\theasciiemail\fi>}
\immediate\write\gtoutfile{\noexpand\\}
\immediate\write\gtoutfile{Authors: \ifx\theasciiauthors\relax
\theauthors\else\theasciiauthors\fi}
{\def\\{ }\immediate\write\gtoutfile{Title: \ifx\theasciititle\relax
\thetitle\else\theasciititle\fi}}
\immediate\write\gtoutfile{Subj-class: GT or SG, GR etc}
\immediate\write\gtoutfile{MSC-class: \theprimaryclass\ifx\thesecondaryclass\relax\else, \thesecondaryclass\fi}
\immediate\write\gtoutfile{Journal-ref: Algebr. Geom. Topol. \thevolumenumber\s
(\thevolumeyear) \startpage-\finishpage}
\immediate\write\gtoutfile{Comments: Published by Algebraic and
Geometric Topology at}
\immediate\write\gtoutfile{\s\s\s  http://www.maths.warwick.ac.uk/agt/AGTVol\thevolumenumber/agt-\thevolumenumber-\thepapernumber.abs.html}
\immediate\write\gtoutfile{\noexpand\\}
\immediate\write\gtoutfile{}
\ifx\theasciiabstract\relax
\immediate\write\gtoutfile{\theabstract}\else
\immediate\write\gtoutfile{\theasciiabstract}\fi
\immediate\write\gtoutfile{}
\immediate\write\gtoutfile{\noexpand\\}
\immediate\write\gtoutfile{}
\immediate\closeout\gtoutfile}}  %%% end of definition of \makeheadfile

\def\maketitlepage{\makeagttitle\makeheadfile}
\let\makeshorttitle\maketitlepage
\let\maketitle\maketitlepage

\def\ifplaintex{\expandafter\ifx\csname documentclass\endcsname\relax}

\def\gtp{{\mathsurround=0pt\it $\cal G\mskip-2mu$eometry \&\ 
$\cal T\!\!$opology $\cal P\!$ublications}}  % GT publications

\def\recd{{\small Received:\qua\receiveddate\ifx\reviseddate\relax
\else\qquad Revised:\qua\reviseddate\fi\par}} 

\def\lognumber#1{\def\thelognumber{#1}}
\def\volumenumber#1{\def\thevolumenumber{#1}}
\def\volumeyear#1{\def\thevolumeyear{#1}}
\def\papernumber#1{\def\thepapernumber{#1}}
\def\pagenumbers#1#2{\def\startpage{#1}\def\finishpage{#2}}
\def\published#1{\def\publishdate{#1}}

\def\received#1{\def\receiveddate{#1}}
\def\revised#1{\def\reviseddate{#1}}
\def\accepted#1{\def\accepteddate{#1}}

\def\asciiaddress#1{\def\theasciiaddress{#1}}

\long\def\asciiabstract#1{\long\def\theasciiabstract{#1}}

\let\\\par\let\thelognumber\relax\let\thevolumenumber\relax
\let\thepapernumber\relax\let\thevolumeyear\relax\let\startpage\relax
\let\finishpage\relax\let\publishdate\relax\let\receiveddate\relax
\let\reviseddate\relax\let\accepteddate\relax\let\theasciititle\relax
\let\theasciiauthors\relax\let\theasciiaddress\relax
\let\theasciiabstract\relax

\let\theasciiemail\relax

\ifplaintex
\font\logobig=cmssbx10 scaled 3836
\font\logomed=cmssbx10 scaled 2557
\else
\font\logobig=cmssbx10 scaled 4200
\font\logomed=cmssbx10 scaled 2800
\fi

\long\def\makeagttitle{   %%% start of definition of \makeagttitle
\count0=\startpage
\agt\hfill      %   Journal title (top left) 
\hbox to 45truept{\vbox to 0pt{\vglue -13truept{\logomed A\kern -.37em{\logobig 
T}\kern -.38em G}\vss}\hss}
\break
{\small Volume \thevolumenumber\ (\thevolumeyear)
\startpage--\finishpage\nl
Published: \publishdate}

\vglue .25truein

{\parskip=0pt\leftskip 0pt plus
1fil\def\\{\par\smallskip}{\Large\bf\thetitle}\par\medskip} \vglue
0.05truein

{\parskip=0pt\leftskip 0pt plus 1fil\def\\{\par}{\sc\theauthors}
\par\medskip}%
 
\vglue 0.03truein 

{\small\leftskip 25truept\rightskip 25truept{\bf Abstract}\stdspace\theabstract

{\bf AMS Classification}\stdspace\theprimaryclass
\ifx\thesecondaryclass\relax\else; \thesecondaryclass\fi\par
{\bf Keywords}\stdspace \thekeywords\par}\vglue 7truept

}   %%%% end of definition of \makeagttitle

\ifplaintex
\font\phead=cmsl9 scaled 950
\font\pnum=cmbx10 scaled 913
\font\pfoot=cmsl9 scaled 950
\headline{\vbox to 0pt{\vskip -4.5mm\line{\small\phead\ifnum
\count0=\startpage ISSN 1472-2739 (on-line) 1472-2747 (printed)
\hfill {\pnum\folio}\else\ifodd\count0\def\\{ }% 
\ifx\theshorttitle\relax\thetitle\else\theshorttitle\fi\hfill{\pnum\folio}
\else\def\\{ and }{\pnum\folio}\hfill\ifx\theshortauthors\relax\theauthors
\else\theshortauthors\fi\fi\fi}\vss}}
\footline{\vbox to 0pt{\vglue 0mm\line{\small\pfoot\ifnum\count0=\startpage
\copyright\ \gtp\hfill\else
\agt, Volume \thevolumenumber\ (\thevolumeyear)\hfill\fi}\vss}}
\else
\font=cmsl9 scaled 1050
\font\lnum=cmbx10 
\font=cmsl9 scaled 1050
\makeatletter
\def\@oddhead{{\small\lhead\ifnum\count0=\startpage ISSN 1472-2739 
(on-line) 1472-2747 (printed)\hfill {\lnum\number\count0}\else\ifodd\count0
\def\\{ }\ifx\theshorttitle\relax \thetitle \else\theshorttitle\fi\hfill
{\lnum\number\count0}\else\def\\{ and }{\lnum\number\count0}
\hfill\ifx\theshortauthors\relax 
\theauthors\else\theshortauthors\fi\fi\fi}}\def\@evenhead{\@oddhead}
\def\@oddfoot{\small\lfoot\ifnum\count0=\startpage\copyright\ \gtp\hfill\else
\agt, Volume \thevolumenumber\ (\thevolumeyear)\hfill\fi}
\def\@evenfoot{\@oddfoot}
\makeatother
\fi
\let\maketitlepage\makeagttitle
\let\makeshorttitle\maketitlepage
\let\maketitle\maketitlepage

\newwrite\gtoutfile
\long\gdef\makeheadfile{  %%% start of definition of \makeheadfile
{\def\\{, }\def\s{ }
\immediate\openout\gtoutfile head.xxx
\immediate\write\gtoutfile{To: math@arxiv.org}
\immediate\write\gtoutfile{Subject: put OR rep NNNNN:ppppp}
\immediate\write\gtoutfile{--text follows this line--}
\immediate\write\gtoutfile{Proxy-for: \ifx\theasciiauthors\relax
\theauthors\else\theasciiauthors\fi\s<\ifx\theasciiemail\relax\theemail\else\theasciiemail\fi>}
\immediate\write\gtoutfile{\noexpand\\}
\immediate\write\gtoutfile{Authors: \ifx\theasciiauthors\relax
\theauthors\else\theasciiauthors\fi}
{\def\\{ }\immediate\write\gtoutfile{Title: \ifx\theasciititle\relax
\thetitle\else\theasciititle\fi}}
\immediate\write\gtoutfile{Subj-class: GT or SG, GR etc}
\immediate\write\gtoutfile{MSC-class: \theprimaryclass\ifx\thesecondaryclass\relax\else, \thesecondaryclass\fi}
\immediate\write\gtoutfile{Journal-ref: Algebr. Geom. Topol. \thevolumenumber\s
(\thevolumeyear) \startpage-\finishpage}
\immediate\write\gtoutfile{Comments: Published by Algebraic and
Geometric Topology at}
\immediate\write\gtoutfile{\s\s\s  http://www.maths.warwick.ac.uk/agt/AGTVol\thevolumenumber/agt-\thevolumenumber-\thepapernumber.abs.html}
\immediate\write\gtoutfile{\noexpand\\}
\immediate\write\gtoutfile{}
\ifx\theasciiabstract\relax
\immediate\write\gtoutfile{\theabstract}\else
\immediate\write\gtoutfile{\theasciiabstract}\fi
\immediate\write\gtoutfile{}
\immediate\write\gtoutfile{\noexpand\\}
\immediate\write\gtoutfile{}
\immediate\closeout\gtoutfile}}  %%% end of definition of \makeheadfile

\def\maketitlepage{\makeagttitle\makeheadfile}
\let\makeshorttitle\maketitlepage
\let\maketitle\maketitlepage

\def\ifplaintex{\expandafter\ifx\csname documentclass\endcsname\relax}

\def\gtp{{\mathsurround=0pt\it $\cal G\mskip-2mu$eometry \&\ 
$\cal T\!\!$opology $\cal P\!$ublications}}  % GT publications

\def\recd{{\small Received:\qua\receiveddate\ifx\reviseddate\relax
\else\qquad Revised:\qua\reviseddate\fi\par}} 

\def\lognumber#1{\def\thelognumber{#1}}
\def\volumenumber#1{\def\thevolumenumber{#1}}
\def\volumeyear#1{\def\thevolumeyear{#1}}
\def\papernumber#1{\def\thepapernumber{#1}}
\def\pagenumbers#1#2{\def\startpage{#1}\def\finishpage{#2}}
\def\published#1{\def\publishdate{#1}}

\def\received#1{\def\receiveddate{#1}}
\def\revised#1{\def\reviseddate{#1}}
\def\accepted#1{\def\accepteddate{#1}}

\def\asciiaddress#1{\def\theasciiaddress{#1}}

\long\def\asciiabstract#1{\long\def\theasciiabstract{#1}}

\let\\\par\let\thelognumber\relax\let\thevolumenumber\relax
\let\thepapernumber\relax\let\thevolumeyear\relax\let\startpage\relax
\let\finishpage\relax\let\publishdate\relax\let\receiveddate\relax
\let\reviseddate\relax\let\accepteddate\relax\let\theasciititle\relax
\let\theasciiauthors\relax\let\theasciiaddress\relax
\let\theasciiabstract\relax

\let\theasciiemail\relax

\ifplaintex
\font\logobig=cmssbx10 scaled 3836
\font\logomed=cmssbx10 scaled 2557
\else
\font\logobig=cmssbx10 scaled 4200
\font\logomed=cmssbx10 scaled 2800
\fi

\long\def\makeagttitle{   %%% start of definition of \makeagttitle
\count0=\startpage
\agt\hfill      %   Journal title (top left) 
\hbox to 45truept{\vbox to 0pt{\vglue -13truept{\logomed A\kern -.37em{\logobig 
T}\kern -.38em G}\vss}\hss}
\break
{\small Volume \thevolumenumber\ (\thevolumeyear)
\startpage--\finishpage\nl
Published: \publishdate}

\vglue .25truein

{\parskip=0pt\leftskip 0pt plus
1fil\def\\{\par\smallskip}{\Large\bf\thetitle}\par\medskip} \vglue
0.05truein

{\parskip=0pt\leftskip 0pt plus 1fil\def\\{\par}{\sc\theauthors}
\par\medskip}%
 
\vglue 0.03truein 

{\small\leftskip 25truept\rightskip 25truept{\bf Abstract}\stdspace\theabstract

{\bf AMS Classification}\stdspace\theprimaryclass
\ifx\thesecondaryclass\relax\else; \thesecondaryclass\fi\par
{\bf Keywords}\stdspace \thekeywords\par}\vglue 7truept

}   %%%% end of definition of \makeagttitle

\ifplaintex
\font\phead=cmsl9 scaled 950
\font\pnum=cmbx10 scaled 913
\font\pfoot=cmsl9 scaled 950
\headline{\vbox to 0pt{\vskip -4.5mm\line{\small\phead\ifnum
\count0=\startpage ISSN 1472-2739 (on-line) 1472-2747 (printed)
\hfill {\pnum\folio}\else\ifodd\count0\def\\{ }% 
\ifx\theshorttitle\relax\thetitle\else\theshorttitle\fi\hfill{\pnum\folio}
\else\def\\{ and }{\pnum\folio}\hfill\ifx\theshortauthors\relax\theauthors
\else\theshortauthors\fi\fi\fi}\vss}}
\footline{\vbox to 0pt{\vglue 0mm\line{\small\pfoot\ifnum\count0=\startpage
\copyright\ \gtp\hfill\else
\agt, Volume \thevolumenumber\ (\thevolumeyear)\hfill\fi}\vss}}
\else
\font=cmsl9 scaled 1050
\font\lnum=cmbx10 
\font=cmsl9 scaled 1050
\makeatletter
\def\@oddhead{{\small\lhead\ifnum\count0=\startpage ISSN 1472-2739 
(on-line) 1472-2747 (printed)\hfill {\lnum\number\count0}\else\ifodd\count0
\def\\{ }\ifx\theshorttitle\relax \thetitle \else\theshorttitle\fi\hfill
{\lnum\number\count0}\else\def\\{ and }{\lnum\number\count0}
\hfill\ifx\theshortauthors\relax 
\theauthors\else\theshortauthors\fi\fi\fi}}\def\@evenhead{\@oddhead}
\def\@oddfoot{\small\lfoot\ifnum\count0=\startpage\copyright\ \gtp\hfill\else
\agt, Volume \thevolumenumber\ (\thevolumeyear)\hfill\fi}
\def\@evenfoot{\@oddfoot}
\makeatother
\fi
\let\maketitlepage\makeagttitle
\let\makeshorttitle\maketitlepage
\let\maketitle\maketitlepage

\newwrite\gtoutfile
\long\gdef\makeheadfile{  %%% start of definition of \makeheadfile
{\def\\{, }\def\s{ }
\immediate\openout\gtoutfile head.xxx
\immediate\write\gtoutfile{To: math@arxiv.org}
\immediate\write\gtoutfile{Subject: put OR rep NNNNN:ppppp}
\immediate\write\gtoutfile{--text follows this line--}
\immediate\write\gtoutfile{Proxy-for: \ifx\theasciiauthors\relax
\theauthors\else\theasciiauthors\fi\s<\ifx\theasciiemail\relax\theemail\else\theasciiemail\fi>}
\immediate\write\gtoutfile{\noexpand\\}
\immediate\write\gtoutfile{Authors: \ifx\theasciiauthors\relax
\theauthors\else\theasciiauthors\fi}
{\def\\{ }\immediate\write\gtoutfile{Title: \ifx\theasciititle\relax
\thetitle\else\theasciititle\fi}}
\immediate\write\gtoutfile{Subj-class: GT or SG, GR etc}
\immediate\write\gtoutfile{MSC-class: \theprimaryclass\ifx\thesecondaryclass\relax\else, \thesecondaryclass\fi}
\immediate\write\gtoutfile{Journal-ref: Algebr. Geom. Topol. \thevolumenumber\s
(\thevolumeyear) \startpage-\finishpage}
\immediate\write\gtoutfile{Comments: Published by Algebraic and
Geometric Topology at}
\immediate\write\gtoutfile{\s\s\s  http://www.maths.warwick.ac.uk/agt/AGTVol\thevolumenumber/agt-\thevolumenumber-\thepapernumber.abs.html}
\immediate\write\gtoutfile{\noexpand\\}
\immediate\write\gtoutfile{}
\ifx\theasciiabstract\relax
\immediate\write\gtoutfile{\theabstract}\else
\immediate\write\gtoutfile{\theasciiabstract}\fi
\immediate\write\gtoutfile{}
\immediate\write\gtoutfile{\noexpand\\}
\immediate\write\gtoutfile{}
\immediate\closeout\gtoutfile}}  %%% end of definition of \makeheadfile

\def\maketitlepage{\makeagttitle\makeheadfile}
\let\makeshorttitle\maketitlepage
\let\maketitle\maketitlepage

\def\ifplaintex{\expandafter\ifx\csname documentclass\endcsname\relax}

\def\gtp{{\mathsurround=0pt\it $\cal G\mskip-2mu$eometry \&\ 
$\cal T\!\!$opology $\cal P\!$ublications}}  % GT publications

\def\recd{{\small Received:\qua\receiveddate\ifx\reviseddate\relax
\else\qquad Revised:\qua\reviseddate\fi\par}} 

\def\lognumber#1{\def\thelognumber{#1}}
\def\volumenumber#1{\def\thevolumenumber{#1}}
\def\volumeyear#1{\def\thevolumeyear{#1}}
\def\papernumber#1{\def\thepapernumber{#1}}
\def\pagenumbers#1#2{\def\startpage{#1}\def\finishpage{#2}}
\def\published#1{\def\publishdate{#1}}

\def\received#1{\def\receiveddate{#1}}
\def\revised#1{\def\reviseddate{#1}}
\def\accepted#1{\def\accepteddate{#1}}

\def\asciiaddress#1{\def\theasciiaddress{#1}}

\long\def\asciiabstract#1{\long\def\theasciiabstract{#1}}

\let\\\par\let\thelognumber\relax\let\thevolumenumber\relax
\let\thepapernumber\relax\let\thevolumeyear\relax\let\startpage\relax
\let\finishpage\relax\let\publishdate\relax\let\receiveddate\relax
\let\reviseddate\relax\let\accepteddate\relax\let\theasciititle\relax
\let\theasciiauthors\relax\let\theasciiaddress\relax
\let\theasciiabstract\relax

\let\theasciiemail\relax

\ifplaintex
\font\logobig=cmssbx10 scaled 3836
\font\logomed=cmssbx10 scaled 2557
\else
\font\logobig=cmssbx10 scaled 4200
\font\logomed=cmssbx10 scaled 2800
\fi

\long\def\makeagttitle{   %%% start of definition of \makeagttitle
\count0=\startpage
\agt\hfill      %   Journal title (top left) 
\hbox to 45truept{\vbox to 0pt{\vglue -13truept{\logomed A\kern -.37em{\logobig 
T}\kern -.38em G}\vss}\hss}
\break
{\small Volume \thevolumenumber\ (\thevolumeyear)
\startpage--\finishpage\nl
Published: \publishdate}

\vglue .25truein

{\parskip=0pt\leftskip 0pt plus
1fil\def\\{\par\smallskip}{\Large\bf\thetitle}\par\medskip} \vglue
0.05truein

{\parskip=0pt\leftskip 0pt plus 1fil\def\\{\par}{\sc\theauthors}
\par\medskip}%
 
\vglue 0.03truein 

{\small\leftskip 25truept\rightskip 25truept{\bf Abstract}\stdspace\theabstract

{\bf AMS Classification}\stdspace\theprimaryclass
\ifx\thesecondaryclass\relax\else; \thesecondaryclass\fi\par
{\bf Keywords}\stdspace \thekeywords\par}\vglue 7truept

}   %%%% end of definition of \makeagttitle

\ifplaintex
\font\phead=cmsl9 scaled 950
\font\pnum=cmbx10 scaled 913
\font\pfoot=cmsl9 scaled 950
\headline{\vbox to 0pt{\vskip -4.5mm\line{\small\phead\ifnum
\count0=\startpage ISSN 1472-2739 (on-line) 1472-2747 (printed)
\hfill {\pnum\folio}\else\ifodd\count0\def\\{ }% 
\ifx\theshorttitle\relax\thetitle\else\theshorttitle\fi\hfill{\pnum\folio}
\else\def\\{ and }{\pnum\folio}\hfill\ifx\theshortauthors\relax\theauthors
\else\theshortauthors\fi\fi\fi}\vss}}
\footline{\vbox to 0pt{\vglue 0mm\line{\small\pfoot\ifnum\count0=\startpage
\copyright\ \gtp\hfill\else
\agt, Volume \thevolumenumber\ (\thevolumeyear)\hfill\fi}\vss}}
\else
\font=cmsl9 scaled 1050
\font\lnum=cmbx10 
\font=cmsl9 scaled 1050
\makeatletter
\def\@oddhead{{\small\lhead\ifnum\count0=\startpage ISSN 1472-2739 
(on-line) 1472-2747 (printed)\hfill {\lnum\number\count0}\else\ifodd\count0
\def\\{ }\ifx\theshorttitle\relax \thetitle \else\theshorttitle\fi\hfill
{\lnum\number\count0}\else\def\\{ and }{\lnum\number\count0}
\hfill\ifx\theshortauthors\relax 
\theauthors\else\theshortauthors\fi\fi\fi}}\def\@evenhead{\@oddhead}
\def\@oddfoot{\small\lfoot\ifnum\count0=\startpage\copyright\ \gtp\hfill\else
\agt, Volume \thevolumenumber\ (\thevolumeyear)\hfill\fi}
\def\@evenfoot{\@oddfoot}
\makeatother
\fi
\let\maketitlepage\makeagttitle
\let\makeshorttitle\maketitlepage
\let\maketitle\maketitlepage

\newwrite\gtoutfile
\long\gdef\makeheadfile{  %%% start of definition of \makeheadfile
{\def\\{, }\def\s{ }
\immediate\openout\gtoutfile head.xxx
\immediate\write\gtoutfile{To: math@arxiv.org}
\immediate\write\gtoutfile{Subject: put OR rep NNNNN:ppppp}
\immediate\write\gtoutfile{--text follows this line--}
\immediate\write\gtoutfile{Proxy-for: \ifx\theasciiauthors\relax
\theauthors\else\theasciiauthors\fi\s<\ifx\theasciiemail\relax\theemail\else\theasciiemail\fi>}
\immediate\write\gtoutfile{\noexpand\\}
\immediate\write\gtoutfile{Authors: \ifx\theasciiauthors\relax
\theauthors\else\theasciiauthors\fi}
{\def\\{ }\immediate\write\gtoutfile{Title: \ifx\theasciititle\relax
\thetitle\else\theasciititle\fi}}
\immediate\write\gtoutfile{Subj-class: GT or SG, GR etc}
\immediate\write\gtoutfile{MSC-class: \theprimaryclass\ifx\thesecondaryclass\relax\else, \thesecondaryclass\fi}
\immediate\write\gtoutfile{Journal-ref: Algebr. Geom. Topol. \thevolumenumber\s
(\thevolumeyear) \startpage-\finishpage}
\immediate\write\gtoutfile{Comments: Published by Algebraic and
Geometric Topology at}
\immediate\write\gtoutfile{\s\s\s  http://www.maths.warwick.ac.uk/agt/AGTVol\thevolumenumber/agt-\thevolumenumber-\thepapernumber.abs.html}
\immediate\write\gtoutfile{\noexpand\\}
\immediate\write\gtoutfile{}
\ifx\theasciiabstract\relax
\immediate\write\gtoutfile{\theabstract}\else
\immediate\write\gtoutfile{\theasciiabstract}\fi
\immediate\write\gtoutfile{}
\immediate\write\gtoutfile{\noexpand\\}
\immediate\write\gtoutfile{}
\immediate\closeout\gtoutfile}}  %%% end of definition of \makeheadfile

\def\maketitlepage{\makeagttitle\makeheadfile}
\let\makeshorttitle\maketitlepage
\let\maketitle\maketitlepage